\documentclass[10pt]{article}
\usepackage{bm,latexsym,amsmath,amsthm,amssymb}
\usepackage{graphicx,indentfirst}
\date{}
\begin{document}
\title{\bf A diffusive logistic problem with a free boundary in time-periodic environment: favorable habitat or unfavorable habitat\thanks{This work was supported by the Fundamental Research Funds for the Central Universities}}
\author{ Qiaoling Chen, Fengquan Li\thanks{Corresponding author.\newline
\mbox{}\qquad E-mail: fqli@dlut.edu.cn (F. Li); qiaolingf@126.com (Q. Chen); mathwangfeng@126.com (F. Wang) },
Feng Wang\\
\small School of Mathematical Sciences, Dalian University of Technology,
Dalian 116024, PR China}
\date{}
\maketitle \baselineskip 5pt
\begin{center}
\begin{minipage}{130mm}
{{\bf Abstract.} We study the diffusive logistic equation with a free boundary in time-periodic environment. To understand the effect of the dispersal rate $d$, the original habitat radius $h_0$, the spreading capability $\mu$ , and the initial density $u_0$ on the dynamics of the problem, we divide the time-periodic habitat into two cases: favorable habitat and unfavorable habitat. By choosing $d, h_0, \mu$ and $u_0$ as variable parameters, we obtain a spreading-vanishing dichotomy and sharp criteria for the spreading and vanishing in time-periodic environment. We introduce the principal eigenvalue $\lambda_1(d, \alpha, \gamma, h(t), T)$ to determine the spreading and vanishing of the invasive species. We prove that if  $\lambda_1(d, \alpha, \gamma, h_0, T)\leq 0$, the spreading must happen; while if $\lambda_1(d, \alpha, \gamma, h_0, T)>0$, the spreading is also possible. Our results show that the species in the favorable habitat can establish itself if the dispersal rate is slow or the occupying habitat is large. In an unfavorable habitat, the species vanishes if the initial density of the species is small, while survive successfully if the initial value is big. Moreover, when spreading occurs, the asymptotic spreading speed of the free boundary is determined.

}

\end{minipage}
\end{center}

\baselineskip=15pt

\section{Introduction}

In this paper, we study the behavior of the solution $(u(t, r), h(t)) (r=|x|, x\in R^N, N\geq2)$, to the following diffusive logistic problem with a free boundary in the time-periodic heterogeneous environment
\begin{align*}
\left\{\begin{array}{l}
u_t-d\Delta u=u(\alpha(t, r)-\gamma(t, r)-\beta(t, r)u), \quad t>0,\quad 0<r<h(t),\\[5pt]
u_{r}(t, 0)=0, u(t, h(t))=0,\quad t>0,\\[5pt]
h'(t)=-\mu u_{r}(t, h(t)), \quad t>0,\\[5pt]
h(0)=h_0, ~u(0, r)=u_0(r), \quad 0\leq r\leq h_{0},
\end{array}\right.
 \tag{1.1}
\end{align*}
where $\triangle u=u_{rr}+\frac{N-1}{r}u_r$; $r=h(t)$ is the free boundary to be determined; $d, h_0$ and $\mu$ are positive constants; the initial function $u_0(r)$ satisfies
\begin{align*}
u_0\in C^2([0, h_0]),\quad u_0'(0)=u_0(h_0)=0,\quad u_0(r)>0 \quad \mbox{in}~ (0, h_0);
 \tag{1.2}
\end{align*}
$\alpha(t, r), \gamma(t, r)$ and $\beta(t, r)$ account for the birth rate, death rate and crowding strength of the species at $r$, respectively, which satisfy the following conditions
\begin{align*}
\left\{\begin{array}{l}
(i)\quad \alpha, \gamma\in C^{\frac{\nu_0}{2}, 1}(R\times[0, \infty)),
 \beta\in C^{\frac{\nu_0}{2}, \nu_0}(R\times[0, \infty))\quad\mbox{for some}~ \nu_0\in(0,1)\\[5pt]
 \qquad \mbox{and are T-periodic in}~t~ \mbox{for some}~T>0;\\[5pt]
(ii)\quad \mbox{there are positive H\"{o}lder continuous and T-periodic functions}~\alpha_i, \gamma_i, \beta_i (i=1, 2)\\[5pt]
 \qquad\mbox{such that}~\alpha_1(t)\leq \alpha(t, r)\leq \alpha_2(t), \gamma_1(t)\leq\gamma(t, r)\leq\gamma_2(t), \\[5pt]
\qquad \beta_1(t)\leq\beta(t, r)\leq\beta_2(t), \quad \forall t\in[0, T], r\in[0, \infty).
\end{array}\right.
 \tag{1.3}
\end{align*}
The equation $h'(t)=-\mu u_{r}(t, h(t))$ in $(1.1)$ is a special case of the well-known Stefan condition, which has been used in the modeling of a number of applied problems \cite{cf00,jc84, lir71}.
We remark that similar free boundary conditions to $(1.3)$ have been used in ecological models over bounded spatial domains in several earlier papers, for example, \cite{l07,myy85,myy86,myy87}.

Recently, Du and Guo \cite{dgp13} considered the following diffusive logistic model with a free boundary in time-periodic environment
\begin{align*}
\left\{\begin{array}{l}
u_t-d\Delta u=u(\alpha(t, r)-\beta(t, r)u), \quad t>0,\quad 0<r<h(t),\\[5pt]
u_{r}(t, 0)=0, u(t, h(t))=0,\quad t>0,\\[5pt]
h'(t)=-\mu u_{r}(t, h(t)), \quad t>0,\\[5pt]
h(0)=h_0, ~u(0, r)=u_0(r), \quad 0\leq r\leq h_{0},
\end{array}\right.
 \tag{1.4}
\end{align*}
where $u(t, |x|)$ represents the population density of an invasive species; $r=h(t)$ denotes the spreading front; $d, h_0$ and $\mu$ are positive constants; $\alpha(t, r)$ and $\beta(t, r)$ are positive T-periodic H\"{o}lder continuous functions, and they satisfy $k_1\leq\alpha(t, r)\leq k_2, k_1\leq\beta(t, r)\leq k_2, \forall t\in[0, T], r\in[0, \infty),$ where $k_1$ and $k_2$ are positive constants. They showed the spreading-vanishing dichotomy in time-periodic environment, and also determined the spreading speed.

In the special case that $\alpha, \gamma$ and $\beta$ are independent of time $t$, problem $(1.1)$ was studied more recently in \cite{llz14}. They showed that, if the diffusion is slow or the occupying habitat is large, the species can establish itself in the favorable habitat, while the species will die out if the initial value of the species is small in an unfavorable habitat. However, big initial number of the species is benefit for the species to survive.

We remark that there are some related recent research about diffusive logistic problem with a free boundary in the homogeneous or heterogeneous environment. In particular, Du and Lin \cite{dl10} are the first ones to study the spreading-vanishing dichotomy of species in the homogeneous environment of dimension one, which has been extended in \cite{dg11} to the situation of higher dimensional space in a radially symmetric case. In \cite{ky11}, some of the results of \cite{dl10} were extended to the case that the solution satisfies a Dirichlet boundary condition at $x=0$ and a free boundary at $x=h(t)$, covering monostable and bistable nonlinearities. Other theoretical advances can also be seen in \cite{dxl13,dmz14,dbl13,pz13,w142,bdk12,zx14}.
Moreover, some Lotka-Vottera competitive type problems with a free boundary were introduced in \cite{wz15,gw12,wz14,dl13}. Other studies of Lotka-Vottera prey-predator problem with a free boundary can be found in \cite{l07,w14,wz13,zw14}.

But the model in \cite{dgp13} do not exactly describe the survival of the species in real environment, for example, some part of the habitat has been polluted or destroyed. To describe the feature of environment, as in \cite{abyn08,lou08}, we say that $r$ is a favorable site if the birth rate $\alpha(t, r)$ is greater than the local death rate $\gamma(t, r)$ over the time interval $[0, T]$, that is, $\int_{0}^{T}(\alpha(t, r)-\gamma(t, r))dt>0$. An unfavorable site is defined by reversing the inequality. Define the favorable set and unfavorable set in $B_{R}$ (a ball with radius $R$) as
\begin{align*}
F_{R}^{+}=\left\{r\in[0, R], \int_{0}^{T}(\alpha(t, r)-\gamma(t, r))dt>0\right\}
\end{align*}
and
\begin{align*}
F_{R}^{-}=\left\{r\in[0, R], \int_{0}^{T}(\alpha(t, r)-\gamma(t, r))dt<0\right\},
\end{align*}
respectively. The habitat $B_{R}$ is characterized as favorable (resp. unfavorable) if the average of the birth rate $\frac{1}{|[0, T]\times B_{R}|}\int_{0}^{T}\int_{B_{R}}\alpha(t, r)dxdt$ is greater than (resp. less than) the average of the death rate $\frac{1}{|[0, T]\times B_{R}|}\int_{0}^{T}\int_{B_{R}}\gamma(t, r)dxdt$.

The aim of this paper is to study the dynamics of problem $(1.1)$ in the time-periodic case, a situation that more closely reflects the periodic variation of the natural environment, such as daily or seasonal changes. Further, we will consider our problem both in the favorable habitat and unfavorable habitat. To best of our knowledge, the present paper seems to be the first attempt to consider the unfavorable habitat with time-periodic environment in the moving domain problem. It should be pointed out here that the arguments developed in the previous works \cite{dgp13} do not work in the situation of unfavorable habitat with time-periodic environment, since now $\alpha(t, r)-\gamma(t, r)$ admitted to change sign in $(t, r)\in[0, T]\times[0, \infty)$. Instead of $h_0$ which is only used in \cite{dgp13}, here we choose the parameters $d, h_0, u_0$ and $\mu$ as the varying parameters to study problem $(1.1)$ both in the favorable habitat and an unfavorable habitat. We derive some sufficient conditions to ensure that spreading and vanishing occur, which yield the spreading-vanishing dichotomy, and sharp criteria governing spreading and vanishing both in the favorable habitat and an unfavorable habitat. Furthermore, we demonstrate that the species will spreading in a favorable habitat if the dispersal rate is slow or the initial occupying habitat radius is large. In an unfavorable habitat, the species will vanish if the initial density of the species is low. However, the species can also spread in an unfavorable habitat, if the initial number of the species is large.
Therefore, we derive that the diffusion, the initial value and the original habitat play a significant role in determining the spreading and vanishing to problem $(1.1)$. Finally, we also extend the asymptotic spreading speed of a free boundary in favorable time-periodic environment, when spreading of the species happens, to unfavorable time-periodic environment.

The rest of our paper is arranged as follows. In section 2, we exhibit some fundamental
results, including the global existence and uniqueness of the solution to problem $(1.1)$ and the comparison principle in the moving domain. An eigenvalue problem is given in section 3. In section 4, a spreading-vanishing dichotomy is proved. Section 5 is devoted to show the sharp criteria governing spreading and vanishing. In section 6, we investigate the asymptotic spreading speed of the free boundary when spreading occurs.

\section{Preliminaries}

In this section, we first give the existence and uniqueness of a global solution to problem $(1.1)$, then a comparison principle which can be used to estimate both $u(t, r)$ and the free boundary $r=h(t)$ is given.

\noindent\textbf{Theorem 2.1.}
For any given $u_0$ satisfying $(1.2)$ and $\alpha, \gamma, \beta$ satisfying $(1.3)$, problem $(1.1)$ admits a unique solution $(u, h)$ defined for all $t>0$. Moreover $h\in C^1([0, \infty)), u\in C^{1, 2}(D)$, with $D=\{(t, r): t>0, 0\leq r\leq h(t)\}$, and $h'(t)>0$ for $t>0, u(t, r)>0$ for $t>0$ and $0\leq r<h(t)$. 

\noindent\textbf{Proof.}
Denote the nonlinear term in $(1.1)$ by a function $f(t, r, u)=u(\alpha(t, r)-\gamma(t, r)-\beta(t, r)u)$. Under the assumption $(1.3)$, it is easy to check that $(i)~ f(t, r, u)$ is H\"{o}lder continuous for $(t, r)\in [0, \infty)\times[0, \infty)$; $(ii)~ f(t, r, u)$ is locally Lipschitz in $u$ uniformly for $(t, r)\in[0, T]\times[0, \infty)$; $(iii)~f(t, r, u)\\
\leq u\max_{[0, T]}\alpha_2(t), \forall (t, r, u)\in[0, T]\times[0, \infty)\times[0, \infty)$.
Therefore, the function $f(t, r, u)$ satisfies conditions in Theorem 3.1 in \cite{dgp13}, then we can prove the local existence and uniqueness result by the contraction mapping theorem, and use suitable estimates to show that the solution is defined for all $t>0$. The proof is a simple modification of Theorem 4.1 in \cite{dg11} and Theorem 3.1 in \cite{dgp13}, so we omit it.
\quad$\Box$

For later applications, we need a comparison principle as follows.

\noindent\textbf{Lemma 2.1.}
Suppose that $T_0\in(0, \infty), \bar{h}\in C^1([0, T_0]), \bar{u}\in C^{1, 2}(D_{T_0})$
with $D_{T_0}=\{(t, r)\in R^2:0\leq t\leq T_0, 0\leq R\leq \bar{h}(t)\}$, and
\begin{align*}
\left\{\begin{array}{l}
\bar{u}_t-d\Delta \bar{u}\geq\bar{u}(\alpha(t, r)-\gamma(t, r)-\beta(t, r)\bar{u}), \quad t>0,\quad 0<r<\bar{h}(t),\\[5pt]
\bar{u}_{r}(t, 0)\leq0, \bar{u}(t, \bar{h}(t))=0,\quad t>0,\\[5pt]
\bar{h}'(t)\geq-\mu \bar{u}_{r}(t, \bar{h}(t)), \quad t>0.\\[5pt]
\end{array}\right.
\end{align*}
If $h_0\leq \bar{h}(0)$ and $u_0(r)\leq \bar{u}(0, r)$ in $[0, h_0]$, then the solution $(u, h)$ satisfies
\begin{align*}
h(t)\leq \bar{h}(t),~ u(t, r)\leq\bar{u}(t, r)\quad \mbox{for}~t\in[0, T_0]
~ \mbox{and}~r\in[0, h(t)].
\end{align*}

\noindent\textbf{Proof.}
The proof is similar to that of Lemma 3.5 in \cite{dl10} and Lemma 2.6 in \cite{dg11} with obvious modification. So we omit it. \quad$\Box$

\noindent\textbf{Remark 2.1.}
The pair $(\bar{u}, \bar{h})$ in Lemma 2.1 is usually called an upper solution of problem $(1.1)$. We can define a lower solution by reversing all the inequalities in the obvious places. Moreover, one can easily prove an analogue of Lemma 2.1 for lower solutions.

\section{An eigenvalue problem}

In this section, we mainly study an eigenvalue problem and analyze the property
of its principle eigenvalue. These results play an important role in later sections. 

Consider the following eigenvalue problem
\begin{align*}
\left\{\begin{array}{l}
\varphi_t-d\Delta \varphi=\left(\alpha(t, |x|)-\gamma(t, |x|)\right)\varphi+\lambda\varphi, \quad \mbox{in}~[0, T]\times B_R,\\[5pt]
\varphi=0,\quad \mbox{on}~[0, T]\times \partial B_R,\\[5pt]
\varphi(0, x)=\varphi(T, x) \quad \mbox{in}~B_R.
\end{array}\right.
 \tag{3.1}
\end{align*}
It is well known \cite{cc03, hess91} that $(3.1)$ possesses a unique principal eigenvalue $\lambda_1=\lambda_1(d, \alpha-\gamma, R, T)$, which corresponds to a positive eigenfunction $\varphi\in C^{1, 2}([0,T]\times \bar{B}_R)$. Moreover, $\varphi(t, x)$ is radially symmetric in $x$ for all t, this fact is a consequence of the moving-plane argument in \cite{dh94}.

In what follows, we present some further properties of $\lambda_1=\lambda_1(d, \alpha- \gamma, R, T)$. We now discuss the dependence of $\lambda_1$ on $d$ for fixed $R$.

\noindent\textbf{Theorem 3.1.}
Let $\alpha(t, r)$ and $\gamma(t, r)$ be functions satisfy $(1.3)$. Then \\
$(i)~\lambda_1(\cdot, \alpha-\gamma, R, T)\rightarrow -\max_{\bar{B}_R}\frac{1}{T}\int_{0}^{T}(\alpha(t, |x|)-\gamma(t, |x|))dt$ as $d\rightarrow 0$;\\
$(ii)~\liminf_{d\rightarrow +\infty}\lambda_1(d, \alpha-\gamma, R, T)=+\infty$.

\noindent\textbf{Proof.}
$(i)$ Rewrite $(3.1)$ in the form
\begin{align*}
\left\{\begin{array}{l}
\varphi_t-d\Delta \varphi=(\alpha(t, |x|)-\gamma(t, |x|)+\lambda_0)\varphi+(\lambda-\lambda_0)\varphi, \quad \mbox{in}~[0, T]\times B_R,\\[5pt]
\varphi=0,\quad \mbox{on}~[0, T]\times \partial B_R,\\[5pt]
\varphi(0, x)=\varphi(T, x) \quad \mbox{in}~B_R,
\end{array}\right.
\end{align*}
where $\lambda_0$ is a constant. We may choose $\lambda_0$ such that
\begin{align*}
-\max_{\bar{B}_R}\frac{1}{T}\int_{0}^{T}(\alpha(t, |x|)-\gamma(t, |x|))dt<\lambda_0
<-\frac{1}{|[0, T]\times B_R|}\int_{0}^{T}\int_{B_R}(\alpha(t, |x|)-\gamma(t, |x|))dxdt.
\end{align*}
For any such $\lambda_0$ we have
\begin{align*}
\frac{1}{|[0, T]\times B_R|}\int_{0}^{T}\int_{B_R}(\alpha(t, |x|)-\gamma(t, |x|)+\lambda_0)dxdt<0
\end{align*}
and
\begin{align*}
\max_{\bar{B}_R}\frac{1}{T}\int_{0}^{T}(\alpha(t, |x|)-\gamma(t, |x|)+\lambda_0)dt>0.
\end{align*}
Similar to Example $17.2^{*}$ in \cite{hess91}, we can apply Prop.~17.1 in \cite{hess91} to prove that $\lambda_1(d, \alpha-\gamma, R, T)-\lambda_0<0$ for small enough $d>0$, which implies
\begin{align*}
\limsup_{d\rightarrow 0} \lambda_1(d, \alpha-\gamma, R, T)\leq -\max_{\bar{B}_R}\frac{1}{T}\int_{0}^{T}(\alpha(t, |x|)-\gamma(t, |x|))dt.
\end{align*}

On the other hand, choosing any $\lambda_0<-\max_{\bar{B}_R}\frac{1}{T}\int_{0}^{T}(\alpha(t, |x|)-\gamma(t, |x|))dt$, we have $\frac{1}{T}\int_{0}^{T}(\alpha(t, |x|)-\gamma(t, |x|)+\lambda_0)dt<0$ for all $x\in \bar{B}_R$. Similar to the proof of Prop.~17.3 in \cite{hess91}, we have $\lambda_1(d, \alpha-\gamma, R, T)-\lambda_0>0$ for small enough $d>0$. Above all, we have $\lambda_1(d, \alpha-\gamma, R, T)\rightarrow -\max_{\bar{B}_R}\frac{1}{T}\int_{0}^{T}(\alpha(t, |x|)-\gamma(t, |x|))dt$ as $d\rightarrow 0$.

$(ii)$
Multiplying $(3.1)$ by $\varphi$ and integrating over $[0, T]\times B_R$, we have
\begin{align*}
\lambda_1(d, \alpha-\gamma, R, T)
&=\frac{\int_{0}^{T}\int_{B_R}[d|\nabla \varphi|^2+
(\gamma(t, |x|)-\alpha(t, |x|))\varphi^2]dxdt}{\int_{0}^{T}\int_{B_R}\varphi^2dxdt}\\
&\geq \frac{d\int_{0}^{T}\int_{B_R}
|\nabla \varphi|^2dxdt}{\int_{0}^{T}\int_{B_R}\varphi^2dxdt}
+\min_{[0,T]\times \bar{B}_{R}}(\gamma(t, |x|)-\alpha(t, |x|)).
\end{align*}
Since $\varphi(t,\cdot)\in H_{0}^{1}(B_{R})$, by Poincar\'{e}'s inequality, there is a constant
$C>0$ depending on $N$ and $R$ such that
$\int_{B_R}|\nabla \varphi|^2dx\geq C\int_{B_R}\varphi^2dx$ for all $t\in [0,T]$, and then
$\int_{0}^{T}\int_{B_R}|\nabla \varphi|^2dxdt\geq k\int_{0}^{T}\int_{B_R}\varphi^2dxdt$.
Thus,
\begin{align*}
\lambda_1(d, \alpha-\gamma, R, T)
\geq dC+\min_{[0,T]\times \bar{B}_{R}}(\gamma(t, |x|)-\alpha(t, |x|)),
\end{align*}
which implies that $\liminf_{d\rightarrow +\infty}\lambda_1(d, \alpha-\gamma, R, T)=+\infty$.
\quad$\Box$

The above theorem implies the following result.

\noindent\textbf{Corollary 3.1.}
$(i)$~If $F_{R}^{+}\neq \emptyset$, then there exists a constant $d_*\in(0, +\infty)$ such that $\lambda_1(d, \alpha- \gamma, R, T)\leq 0$ for $0<d\leq d_*$;
$(ii)$~there exists a large constsnt $d^{*}\in(0, +\infty)$ such that $\lambda_1(d, \alpha-\gamma, R, T)>0$ for $d>d^*$. 

Next, we give the dependence of $\lambda_1(d, \alpha-\gamma, R, T)$ on $R$ for fixed $d$. To ensure that $\lambda_1(d, \alpha-\gamma, R, T)=0$ has a positive root $h^{*}$, we need an extra assumption as follows:
there exist $T-$periodic positive H\"{o}lder continuous functions $\eta_{*}(t)$ and $\eta^{*}(t)$ such that
\begin{align*}
(H) \quad \liminf_{x\rightarrow+\infty}~(\alpha(t,r)-\gamma(t,r))=\eta_{*}(t), \quad
\limsup_{x\rightarrow+\infty}~(\alpha(t,r)-\gamma(t,r))=\eta^{*}(t).
\end{align*}

\noindent\textbf{Theorem 3.2.}
Let $\alpha(t, r)$ and $\gamma(t, r)$ be functions satisfying $(1.3)$. Then\\
$(i)~\lambda_1(d, \alpha-\gamma, \cdot, T)$ is a strictly decreasing continuous function in $(0, +\infty)$ for fixed $d, \alpha, \gamma, T$, and $\lambda_1(d, \cdot, R, T)$ is a strictly decreasing continuous function in the sense that, $\lambda_1(d, k_{1}(t,r), R, T)<\lambda_1(d, k_{2}(t,r), R, T)$ if the two T-periodic continuous functions $k_{1}(t,r)$ and $k_{2}(t,r)$ satisfy $k_{1}(t,r)\geq, \not\equiv k_{2}(t,r)$ on $[0, T]\times B_R$;\\
$(ii)~\lambda_1(d, \alpha-\gamma, R, T)\rightarrow+\infty$ as $R\rightarrow 0$;\\
$(iii)~\lim_{R\rightarrow\infty}\lambda_1(d, \alpha-\gamma, R, T)<0$ under the assumption $(H)$.

\noindent\textbf{Proof.}
$(i)$ The continuity of $\lambda_1(d, \alpha-\gamma, \cdot, T)$ with fixed $d, \alpha, \gamma, T$ can be obtained by using a simple re-scaling argument of the spatial variable $r$, which also gives the monotonicity of $\lambda_1(d, \alpha-\gamma, \cdot, T)$. Further, it follows from \cite{cc03, hess91} that the continuity and monotonicity of $\lambda_1(d, \cdot, R, T)$.

$(ii)$ From the condition $(1.3)$ and the monotonicity of $\lambda_1(d, k, R, T)$ in $k$, we have
\begin{align*}
\alpha(t, r)-\gamma(t, r)\leq\alpha_2(t)\leq\max_{[0, T]}\alpha_2(t),
\end{align*}
which implies
\begin{align*}
\lambda_1(d, \max_{[0, T]}\alpha_2(t), R, T)\leq \lambda_1(d, \alpha-\gamma, R, T).
\end{align*}
It is also easy to check that $\lambda_1(d, \max_{[0, T]}\alpha_2(t), R, T)$ is the principal eigenvalue of the elliptic problem
\begin{align*}
\left\{\begin{array}{l}
-d\Delta \varphi=\max_{[0, T]}\alpha_2(t)\varphi+\lambda\varphi, \quad \mbox{in}~B_R,\\[5pt]
\varphi=0,\quad \mbox{on}~\partial B_R.
\end{array}\right.
\end{align*}
It is well-known that
\begin{align*}
\lim_{R\rightarrow 0^{+}}\lambda_1(d, \max_{[0, T]}\alpha_2(t), R, T)=+\infty.
\end{align*}
Thus $\lim_{R\rightarrow 0^{+}}\lambda_1(d, \alpha-\gamma, R, T)=+\infty$.

$(iii)$ From the monotony of $\lambda_{1}(d, k(t, r), R, T)$ in $k$, we have
\begin{align*}
\lambda_{1}(d, \alpha(t,r)-\gamma(t,r), R, T)
\leq \lambda_{1}(d, \min_{[0,T]}(\alpha(t,r)-\gamma(t,r)), R, T)
\end{align*}
It is also easy to see that $\lambda_{1}(d, \min_{[0,T]}(\alpha(t,r)-\gamma(t,r)), R, T)$ is the principal eigenvalues of the elliptic problems
\begin{align*}
\left\{\begin{array}{l}
-d\Delta \varphi=\min_{[0,T]}(\alpha(t,r)-\gamma(t,r))\varphi+\lambda\varphi, \quad \mbox{in}~B_R,\\[5pt]
\varphi=0,\quad \mbox{on}~\partial B_{R}.
\end{array}\right.
\end{align*}
Note that $\min_{[0,T]}(\alpha(t,r)-\gamma(t,r))$ is continuous in $[0, \infty)$
and satisfies
\begin{align*}
0<\min_{[0,T]}\eta_{*}(t)
\leq\liminf_{x\rightarrow \infty}\min_{[0,T]}(\alpha(t,r)-\gamma(t,r))
\leq\limsup_{x\rightarrow \infty}\min_{[0,T]}(\alpha(t,r)-\gamma(t,r))
\leq \max_{[0,T]}\eta^{*}(t).
\end{align*}
From Remark 3.1 in \cite{w142}, we know
\begin{align*}
\lim_{R\rightarrow \infty}\lambda_{1}(d, \min_{[0,T]}(\alpha(t,r)-\gamma(t,r)), R, T)<0,
\end{align*}
and then
\begin{align*}
\lim_{R\rightarrow \infty}\lambda_{1}(d, (\alpha(t,r)-\gamma(t,r)), R, T)<0.
\end{align*}
This completes the proof.\quad$\Box$

The above theorem implies the following result.

\noindent\textbf{Corollary 3.2} There exists a threshold
$h^{*}\in(0,\infty]$ such that $\lambda_{1}(d, \alpha-\gamma, R, T)\leq 0$
for $R\geq h^{*}$ and $\lambda_{1}(d, \alpha-\gamma, R, T)>0$ for $0<R<h^{*}$.
Moreover, $h^{*}\in(0,\infty)$ if the assumption $(H)$ holds. If we replace $R$ in
$(3.1)$ by $h(t)$, then it follows from the strict increasing monotony of $h(t)$
and Theorem $3.2$ that $\lambda_{1}(d, \alpha-\gamma, h(t), T)$ is a strictly
monotone decreasing function of $t$.

\section{The Spreading-vanishing dichotomy}

In this section, we prove the spreading-vanishing dichotomy. Though our approach here mainly follows the lines of \cite{dgp13},
considerable changes in the proofs are needed, since the situation here is more general and difficult.

In order to study the long time behavior of the spreading species, we firstly consider the fixed boundary problem
\begin{align*}
\left\{\begin{array}{l}
w_{t}-d\Delta w=w(\alpha(t,r)-\gamma(t,r)-\beta(t,r)w), \quad t>0,~ 0<r<R,\\[5pt]
w_{r}(t,0)=0, w(t,R)=0, \quad t>0,\\[5pt]
w(0,r)=u_{0}(r), \quad 0\leq r\leq R,
\end{array}\right.
 \tag{4.1}
\end{align*}
where $R$ is a positive constant, and $u_{0}(r)$ satisfies $(1.2)$. The related $T-$periodic problem is
\begin{align*}
\left\{\begin{array}{l}
W_{t}-d\Delta W=W(\alpha(t,r)-\gamma(t,r)-\beta(t,r)W), \quad (t,r)\in[0,T]\times [0,R],\\[5pt]
W_{r}(t,0)=0, W(t,R)=0, \quad t\in[0,T],\\[5pt]
W(0,r)=W(T,r), \quad r\in[0,R].
\end{array}\right.
 \tag{4.2}
\end{align*}

\noindent\textbf{Proposition 4.1.}
$(i)$ if $\lambda_{1}(d, \alpha-\gamma, R, T)<0$, problem (4.2) has a unique $T-$periodic
positive solution $W(t,r)$;\\
$(ii)$ if $\lambda_{1}(d, \alpha-\gamma, R, T)<0$, any positive solution of problem (4.1)
converges to $W(t,r)$ in $C_{0}^{1}(\overline{B_{R}})$ as $t\rightarrow \infty$;\\
$(iii)$ if $\lambda_{1}(d, \alpha-\gamma, R, T)\geq0$, any positive solution of problem (4.1) converges to $0$ in $C_{0}^{1}(\overline{B_{R}})$ as $t\rightarrow \infty$.

\noindent\textbf{Proof.} the proof of $(i)$ and $(ii)$ can be see in \cite{cc03,hess91}.
Taking $p=2$ in Theorem 3.1 in \cite{pw12}, we can easily get $(iii)$.\quad$\Box$

Next we consider the $T-$periodic parabolic logistic problem on the entire space $\mathbb{R}^{N}(N\geq 2)$
\begin{align*}
\left\{\begin{array}{l}
U_{t}-d\Delta U=U(\alpha(t,|x|)-\gamma(t,|x|)-\beta(t,|x|)U), \quad (t,x)\in[0,T]\times \mathbb{R}^{N},\\[5pt]
U(0,|x|)=U(T,|x|), \quad x\in\mathbb{R}^{N}.
\end{array}\right.
 \tag{4.3}
\end{align*}

\noindent\textbf{Proposition 4.2.}
Assume that $(H)$ holds, then\\
$(i)$ problem $(4.3)$ has a unique positive $T-$periodic (radial) solution $U(t,r)$;\\
$(ii)$ the unique positive $T-$periodic solution
of problem $(4.2)$ satisfies $W(t,r)\rightarrow U(t,r)$ as $R\rightarrow +\infty$.

\noindent\textbf{Proof.}
Taking $p=2$ in Theorem 1.3 in \cite{pw12}, we can easily get $(i)$ and $(ii)$.\quad$\Box$

It follows from Theorem 2.1 that $h(t)$ is monotonic increasing and therefore there exists $h_{\infty}\in (0,+\infty]$ such
that $\lim_{t\rightarrow+\infty}h(t)=h_{\infty}$. The spreading-vanishing dichotomy  is a consequence of the following two lemmas.

\noindent\textbf{Lemma 4.1} If $h_{\infty}<\infty$, then $\lim_{t\rightarrow+\infty}\|u(t,\cdot)\|_{C([0,h(t)])}=0$. Furthermore, we have $h_{\infty}\leq h^{*}$ under the extra assumption $(H)$.

\noindent\textbf{Proof.}
Firstly, we use the contradiction argument to prove $\lim_{t\rightarrow+\infty}\|u(t,\cdot)\|_{C([0,h(t)])}=0$.

Support that $\varepsilon_{0}=\limsup_{t\rightarrow\infty}\|u(t, \cdot)\|_{C([0, h(t)])}>0$, then there exists a sequence
$(t_k, r_k)\in (0,\infty)\times [0,h(t)]$ with $t_{k}\rightarrow \infty$ as $k\rightarrow \infty$ such that $u(t_k, x_k)\geq\frac{\varepsilon_{0}}{2}$ for all $k\in N$. Since $0\leq r_k<h_{\infty}<\infty$, passing to a subsequence if necessary,
one has $r_{k}\rightarrow r_{0}$ as $k\rightarrow \infty$. We claim $r_{0}\in[0,h_{\infty})$. Indeed, if $r_{0}=h_{\infty}$,
then $r_k-h(t_k)\rightarrow 0$ as $k\rightarrow\infty$. Furthermore, we have
\begin{align*}
\left|\frac{\varepsilon_{0}}{2(r_k-h(t_k))}\right|
\leq\left|\frac{u(t_k, r_k)}{r_k-h(t_k)}\right|
=\left|\frac{u(t_k, r_k)-u(t_k, h(t_k))}{r_k-h(t_k)}\right|
=|u_{r}(t_k, \bar{r}_k)|\leq C,
\end{align*}
where $\bar{r}_k\in(r_k, h(t_k))$. It is a contradiction since $r_k-h(t_k)\rightarrow 0$ as $k\rightarrow \infty$.

Define
\begin{align*}
u_k(t, r)=u(t+t_k, r) \quad \mbox{for}~(t, x)\in D_k,
\end{align*}
with $D_k:=\{(t, r)\in \mathbb{R}^2: t\in(-t_k, \infty), r\in [0, h(t+t_k)]\}$.

We claim that the unique positive solution $u(t,r)$ to $(1.1)$ is bounded, that is, there exists a positive constant $M$ such that $0<u(t,r)\leq M$ for $t\in (0,+\infty)$ and $0\leq r<h(t)$. In fact,
it follows from the comparison principle that $u(t,r)\leq \bar{u}(t)$ for $t\in[0,\infty)$ and $r\in[0,h(t)]$, where $\bar{u}(t)$
is the solution of the problem
\begin{align*}
\left\{\begin{array}{l}
\frac{d\bar{u}}{dt}=\bar{u}(\max_{[0,T]}\alpha_{2}(t)-\min_{[0,T]}\beta_{1}(t)\bar{u}), \quad t>0,\\[5pt]
\bar{u}(0)=\|u_{0}\|_{L^{\infty}([0,h_{0}])}.
\end{array}\right.
\end{align*}
Thus we have
\begin{align*}
u(t,r)\leq \sup_{t\geq 0}\bar{u}(t)
=\max\left\{\frac{\max_{[0,T]\alpha_{2}(t)}}{\min_{[0,T]\beta_{1}(t)}},\|u_{0}\|_{L^{\infty}([0,h_{0}])}\right\}\triangleq M.
\tag{4.4}
\end{align*}

By the parabolic regularity (see \cite{fri64,liberman96,lsu67}), we have, up to a subsequence if necessary, $u_{k}\rightarrow \hat{u}$ as $k\rightarrow \infty$
with $\hat{u}$ satisfying
\begin{align*}
\hat{u}_{t}-d\Delta \hat{u}=\hat{u}(\alpha(t,r)-\gamma(t,r)-\beta(t,r)\hat{u}), \quad (t,r)\in (-\infty,+\infty)\times(0,h_{\infty}).
\end{align*}
We note that
\begin{align*}
\hat{u}(0,r_{0})=\lim_{k\rightarrow \infty}u_{k}(0,r_{k})=\lim_{k\rightarrow \infty}u_{k}(t_{k},r_{k})\geq\frac{\varepsilon_{0}}{2}.
\end{align*}
By the strong maximum principle, $\hat{u}>0$ in $(-\infty,+\infty)\times[0,h_{\infty})$, and
\begin{align*}
\hat{u}_{t}-d\Delta \hat{u}\geq-\hat{u}\|\alpha(t,r)-\gamma(t,r)-\beta(t,r)\hat{u}\|_{L^{\infty}(\mathbb{R}\times\mathbb{R}_{+})}.
\end{align*}
Furthermore, using the Hopf lemma to the equation of $\hat{u}$ at the point $(0,h_{\infty})$, we obtain $\hat{u}_{x}(0,h_{\infty})<0$.

On the other hand, since $h_{\infty}<\infty$ and the bound of $h^{\prime}(t)$ is independent of $t$ (the proof is a simple modification
of that for Theorem $4.1$ in \cite{dgp13}). We claim that $h^{\prime}(t)\rightarrow 0$
as $t\rightarrow\infty$. In fact, support that there exists a sequence $\{t_{n}\}$ such that $t_{n}\rightarrow \infty$
and $h^{\prime}(t_{n})\rightarrow \varepsilon_{1}$ as $n\rightarrow \infty$ for some $\varepsilon_{1}>0$. We can find $\varepsilon_{2}>0$
small enough such that $h^{\prime}(t)\geq\frac{\varepsilon_{1}}{2}$ for all $t\in(t_{n}-\varepsilon_{2},t_{n}+\varepsilon_{2})$ for
all $n$. Then we obtain
\begin{align*}
h_{\infty}=h_{0}+\int_{0}^{\infty}h^{\prime}(t)dt
\geq h_{0}+\sum_{n=1}^{\infty}\int_{t_{n}-\varepsilon}^{t_{n}+\varepsilon}\frac{\varepsilon_{1}}{2}dt=\infty.
\end{align*}
a contradiction. Hence, we can derive
$u_{x}(t_{k},h(t_{k}))\rightarrow 0$ as $k\rightarrow \infty$ in view of the Stefan condition, i.e. $h'(t)=-\mu u_{r}(t, h(t))$. However,
$u_{x}(t_{k},h(t_{k}))=(u_{k})_{x}(0,h(t_{k}))\rightarrow \hat{u}_{x}(0,h_{\infty})$ as $k\rightarrow \infty$, which produces
a contradiction that $0=\hat{u}_{x}(0,h_{\infty})<0$.
Therefore, $\lim_{t\rightarrow \infty}\|u(t, \cdot)\|_{C([0, h(t)])}=0$.

Next, we claim that $h_{\infty}\leq h^{*}$. It is sufficient to prove $\lambda_{1}(d,\alpha-\gamma,h_{\infty},T)\geq0$ by Corollary 3.2. We argue by contradiction.
By the continuity of $\lambda_{1}(d,\alpha-\gamma,h(t),T)$ in $t$ and $h(t)\rightarrow h_{\infty}$
as $t\rightarrow \infty$, there exists $\tau\gg 1$ such that $\lambda_{1}(d,\alpha-\gamma,h(\tau),T)<0$.
Let $\underline{u}$ be the solution of the problem
\begin{align*}
\left\{\begin{array}{l}
\underline{u}_{t}-d\Delta \underline{u}=\underline{u}(\alpha(t,r)-\gamma(t,r)-\beta(t,r)\underline{u}), \quad t\geq \tau,~ 0<r< h(\tau),\\[5pt]
\underline{u}_{r}(t,0)=0, \underline{u}(t,h(\tau))=0, \quad t\geq \tau,\\[5pt]
\underline{u}(\tau,r)=u(\tau,r), \quad 0\leq r\leq h(\tau).
\end{array}\right.
\end{align*}
Then by the comparison principle (see Remark 2.1), we have $\underline{u}(t,r)\leq u(t,r)$ for $(t,r)\in[\tau,\infty)\times [0,h(\tau)]$. It follows from Proposition 4.1 that $\lim_{n\rightarrow+\infty}\underline{u}(t+nT, r)=\underline{u}^*(t, r)$ uniformly on $[0, h(\tau)]$,
where $\underline{u}^{*}(t,r)$ is the unique positive $T-$periodic solution of the problem
\begin{align*}
\left\{\begin{array}{l}
\underline{u}^{*}_{t}-d\Delta \underline{u}^{*}=\underline{u}^{*}(\alpha(t,r)-\gamma(t,r)-\beta(t,r)\underline{u}^{*}),
\quad (t,r)\in [0,T]\times[0,h(\tau)],\\[5pt]
\underline{u}^{*}_{r}(t,0)=0, \underline{u}^{*}(t,h(\tau))=0, \quad t\in [0, T],\\[5pt]
\underline{u}^{*}(0,r)=u^{*}(T,r), \quad r\in [0,h(\tau)].
\end{array}\right.
\end{align*}
Hence, $\liminf_{n\rightarrow \infty}u(t+nT,r)\geq \underline{u}^{*}(t,r)>0$ in $(0,h(\tau))$.
This contradicts to $\lim_{t\rightarrow+\infty}\|u(t,\cdot)\|_{C([0,h(t)])}\\=0$,
which completes the proof of Lemma 4.1. \quad$\Box$

\noindent\textbf{Lemma 4.2} Assume that $(H)$ holds. If $h_{\infty}=\infty$, then
\begin{align*}
\lim_{n\rightarrow \infty}u(t+nT,r)=U(t,r) \quad\mbox{locally uniformly for}~ (t,r)\in[0,T]\times [0,\infty),
\tag{4.5}
\end{align*}
where $U(t,r)$ is the unique positive $T-$periodic (radial) solution of problem $(4.3)$.

\noindent\textbf{Proof.}
Since $(H)$ holds, we have $\liminf_{r\rightarrow+\infty}(\alpha(t,r)-\gamma(t,r))=\eta_{*}(t)>0$, and then there exists a large positive
constant $r^{*}$ such that $\alpha(t, r)-\gamma(t, r)>0$
and $\lambda_{1}(d,\alpha-\gamma,R,T)<0$ for any $R>r^{*}$ (by
Corollary $3.2$). Thus, we can choose an increasing
sequence of positive numbers $R_{m}$ with $R_{1}>r^{*}$ and $R_{m}\rightarrow+\infty$ as $m\rightarrow+\infty$, such that
$\lambda_{1}(d,\alpha-\gamma,R_{m},T)<0$ for all $m\geq1$. And we note that there exists a $t^{*}>0$ such that $r^{*}=h(t^{*})$.

Let $W_{R_{m}}(t,r)$ be the unique positive $T-$periodic solution to $(4.2)$ with replacing $R$ by $R_{m}$. It follows from proposition $4.2$ that $W_{R_{m}}(t,r)$ converge to $U(t,r)$ as $R_{m}\rightarrow+\infty$. Since $h_{\infty}=\infty$, we can find $T_{m}>0$ such that $h(t)\geq R_{m}$ for all $t\geq T_{m}$. We consider the following problem
\begin{align*}
\left\{\begin{array}{l}
w_{t}-d\Delta w=w(\alpha(t,r)-\gamma(t,r)-\beta(t,r)w), \quad t\geq T_{m},~ 0\leq r\leq R_{m},\\[5pt]
w_{r}(t,0)=0, w(t,R_{m})=0, \quad t\geq T_{m},\\[5pt]
w(T_{m},r)=u(T_{m},r), \quad 0\leq r\leq R_{m}.
\end{array}\right.
 \tag{4.6}
\end{align*}
By the comparison principle, we have $w_{R_{m}}(t,r)\leq u(t,r)$ for $t\geq T_{m}$ and $0\leq r\leq R_{m}$.

Since $\lambda_{1}(d,\alpha-\gamma,R_{m},T)<0$ for all $m\geq1$, it follows from Proposition $4.1$ that all positive solutions
$w_{R_{m}}(t+nT,r)$ of $(4.6)$ converge to $W_{R_{m}}(t,r)$ uniformly for $(t,r)\in[0,T]\times[0,R_{m}]$ as $n\rightarrow+\infty$.
Therefore, $W_{R_{m}}(t,r)\leq \liminf_{n\rightarrow+\infty}u(t+nT,r)$ uniformly for $(t,r)\in[0,T]\times[0,R_{m}]$.
Sending $m\rightarrow \infty$, it follows from Proposition $4.2$ that
\begin{align*}
U(t,r)\leq \liminf_{n\rightarrow+\infty}u(t+nT,r) \quad\mbox{locally uniformly for}~ (t,r)\in[0,T]\times [0,\infty).
\tag{4.7}
\end{align*}

Next, using a squeezing argument similar in spirit to \cite{dm01},
we prove that
\begin{align*}
\limsup_{n\rightarrow+\infty}u(t+nT,r)\leq U(t,r) \quad\mbox{locally uniformly for}~ (t,r)\in[0,T]\times [0,\infty).
\tag{4.8}
\end{align*}
We first consider the following T-periodic boundary blow-up parabolic problems
\begin{align*}
\left\{\begin{array}{l}
V_{t}-d\Delta V=V(\alpha(t+t^{*},r+r^{*})-\gamma(t+t^{*},r+r^{*})-\beta(t+t^{*},r+r^{*})V), ~ (t,r)\in[0,T]\times [0,R_{m}],\\[5pt]
V(t+t^{*},R_{m}+r^{*})=\infty, ~ t\in[0,T],\\[5pt]
V(t^{*},r+r^{*})=V(T+t^{*},r+r^{*}), \quad r\in[0,R_{m}],
\end{array}\right.
 \tag{4.9}
\end{align*}
It follows from Lemma 3.1 in \cite{pw12} that $(4.9)$ has a unique positive $T-$periodic solution $V_{R_{m}}(t+t^{*},r+r^{*})\triangleq V_{R_{m}}^{*}(t,r)$.
Moreover, for any integer
$K\geq \max_{[0,T]\times \overline{B_{R_{m}}}}\frac{\alpha(t+t^{*},r+r^{*})-\gamma(t+t^{*},r+r^{*})}{\beta(t+t^{*},r+r^{*})}$,
we consider the problem
\begin{align*}
\left\{\begin{array}{l}
V_{t}-d\Delta V=V(\alpha(t+t^{*},r+r^{*})-\gamma(t+t^{*},r+r^{*})-\beta(t+t^{*},r+r^{*})V), \quad (t,r)\in[0,T]\times [0,R_{m}],\\[5pt]
V(t+t^{*},R_{m}+r^{*})=K, \quad t\in[0,T],\\[5pt]
V(t^{*},r+r^{*})=V(T+t^{*},r+r^{*}), \quad r\in[0,R_{m}].
\end{array}\right.
 \tag{4.10}
\end{align*}
Since $V=0$ and $V=K$ are the sub-supersolutions to $(4.10)$, then it has one positive solution. By the same analysis as the uniqueness of positive $T-$periodic solution to $(4.2)$, the uniqueness of positive $T-$periodic solution follows,
denote it by $V_{R_{m},K}^{*}(t,r)$. By the comparison principle, it is obvious to check that
$V_{R_{m},K}^{*}(t,r)<V_{R_{m},K+1}^{*}(t,r)$. Thus, $\lim_{K\rightarrow+\infty}V_{R_{m},K}^{*}(t,r)$ exists and is a positive
$T-$periodic solution to $(4.9)$. Since $(4.9)$ has a unique positive $T-$ periodic solution, then
$\lim_{K\rightarrow+\infty}V_{R_{m},K}^{*}(t,r)=V_{R_{m}}^{*}(t,r)$.

Furthermore, let
$\bar{V}_{R_{m}}(t,r)=\bar{V}(t+t^{*},r+r^{*})$ be the positive solution to the following parabolic problem
\begin{align*}
\left\{\begin{array}{l}
\bar{V}_{t}-d\Delta \bar{V}=\bar{V}(\alpha(t+t^{*},r+r^{*})-\gamma(t+t^{*},r+r^{*})-\beta(t+t^{*},r+r^{*})\bar{V}),
\quad t>0,~0\leq r\leq R_{m},\\[5pt]
\bar{V}(t+t^{*}, R_{m}+r^{*})=K, \quad t>0,\\[5pt]
\bar{V}(t^{*},r+r^{*})=k_{m}V_{R_{m},K+1}^{*}(t,r), \quad 0\leq r\leq R_{m},
\end{array}\right.
\end{align*}
where $k_{m}$ is a positive constant and satisfies $k_{m}V_{R_{m},K+1}^{*}(t, R_{m})\geq K$. Since
$\lambda_{1}(d,\alpha-\gamma,R_{m},T)<0$, it follows from Proposition $4.1$ that
$\bar{V}_{R_{m}}(t+nT,r)\rightarrow V_{R_{m},K}^{*}(t,r)$ as $n\rightarrow \infty$ for
$(t,r)\in [0,T]\times [0,R_{m}]$. By the comparison principle,
$u(t+t^{*}, r+r^{*})\leq \bar{V}_{R_{m}}(t,r)$ for $t> 0$ and $0\leq r\leq R_{m}$,
which implies that
$\limsup_{n\rightarrow \infty}u(t+t^{*}+nT,r+r^{*})\leq V_{R_{m},K}^{*}(t,r)$
uniformly for $(t,r)\in [0,T]\times [0,R_{m}]$.
Therefore, from Proposition 4.2, we get
$\limsup_{n\rightarrow \infty}u(t+nT,r)\leq U(t,r)$ locally
uniformly for $(t,r)\in [0,T]\times [0,\infty)$. Clearly, $(4.5)$ is a consequence of $(4.7)$ and $(4.8)$. \quad$\Box$

Combining Lemma $4.1$ and $4.2$, we immediately obtain the following spreading-vanishing dichotomy theorem.

\noindent\textbf{Theorem 4.1.} Assume that $(H)$ holds.
Let $(u(t,r), h(t))$ be a solution of problem $(1.1)$. Then, the following alternative holds:\\
Either~ $(i)$ spreading: $h_{\infty}=\infty$ and $\lim_{n\rightarrow+\infty}u(t+nT,r)=U(t,r)$ locally
uniformly for $(t,r)\in [0,T]\times [0,\infty)$;\\
or~ $(ii)$ vanishing: $h_{\infty}\leq h^{*}<\infty$ and
$\limsup_{t\rightarrow\infty}\|u(t, \cdot)\|_{C([0, h(t)])}=0$,
where $U(t,r)$ is the unique positive $T-$periodic solution to $(4.3)$.

\section{Sharp criteria for spreading and vanishing}

In this section, we aim to use parameters $d$, $h_{0}$,
$\mu$ and $u_0(r)$ to derive sharp criteria for species spreading and vanishing established in Theorem $4.1$.

\subsection{Case $\lambda_{1}(d,\alpha-\gamma,h_{0},T)\leq 0$}

In this subsection, we discuss the case $\lambda_{1}(d,\alpha,\gamma,h_{0},T)\leq 0$.

\noindent\textbf{Lemma 5.1.}
$h_{\infty}=\infty$.

\noindent\textbf{Proof.}
$\lambda_{1}(d,\alpha-\gamma,h_{0},T)\leq 0$ can be considered in two cases. One is $\lambda_{1}(d,\alpha-\gamma,h_{0},T)<0$,
the other is  $\lambda_{1}(d,\alpha-\gamma,h_{0},T)=0$. Note that  $\lambda_{1}(d,\alpha-\gamma,h(t),T)$ is strictly monotone decreasing function of $t$. Therefore,  $\lambda_{1}(d,\alpha-\gamma,h(\tilde{T}),T)<\lambda_{1}(d,\alpha-\gamma, h_{0},T)=0$ for any positive time $\tilde{T}$. Setting $\tilde{T}$ be the initial time, then we go back the first case
$\lambda_{1}(d,\alpha-\gamma,h_{0},T)<0$. Thus, we only show the case of $\lambda_{1}(d,\alpha-\gamma,h_{0},T)<0$. To simplify, we denote $\lambda_{1}=\lambda_{1}(d,\alpha-\gamma,h_{0},T)$, and its corresponding eigenfunction by $\varphi$ with $\|\varphi\|_{L^{\infty}([0,T]\times[0,h_{0}])}=1$.

Now we set
\begin{align*}
\underline{u}(t, r)=\left\{\begin{array}{l}
\varepsilon\varphi(t,r), \quad \mbox{for}~t\geq 0, 0\leq r\leq h_{0},\\[5pt]
0, \quad \mbox{for}~t\geq 0, r>h_{0}.
\end{array}\right.
\end{align*}
Choose $\varepsilon>0$ so small that
\begin{align*}
\varepsilon\leq \min\left\{u_{0}(r),-\frac{\lambda_{1}}{\max_{[0,T]}\beta_{2}(t)}\right\},
\quad \mbox{for}\quad r\in[0, h_{0}].
\end{align*}
Then by a direct calculation, we obtain
\begin{align*}
\left\{\begin{array}{l}
\underline{u}_t-d\Delta\underline{u}-\underline{u}(\alpha(t,r)-\gamma(t,r)-\beta(t,r)\underline{u})
=\varepsilon\varphi(\lambda_1+\beta(t,r)\varepsilon\varphi)\leq0, \quad t>0,\quad 0<r<h_{0},\\[5pt]
\underline{u}_{r}(t, 0)=0, \underline{u}(t, h_{0})=0, \quad t>0,\\[5pt]
0=h_{0}^{\prime}\leq -\mu\underline{u}_{r}(0, h_{0})=-\mu\varepsilon\varphi^{\prime}(h_{0}),
\quad t>0,\\[5pt]
\underline{u}(0, r)\leq u_{0}(r), \quad 0\leq r\leq h_{0}.
\end{array}\right.
\end{align*}
By the comparison principle, we have
\begin{align*}
\underline{u}(t, r)\leq u(t, r)
\quad \mbox{for}~[0, \infty)\times[0, h_{0}].
\end{align*}
It follows that
\begin{align*}
\liminf_{t\rightarrow\infty}\|u(t, \cdot)\|_{C([0, h(t)])}\geq\epsilon\varphi_1(t,0)>0.
\end{align*}
According to Lemma 4.1, we see that $h_{\infty}=\infty$. \quad$\Box$

Assume that there is a favorable site in $B_{h_0}$ and the diffusion is slow ($0<d\leq d_{*}$), by applying Corollary 3.1, we have $\lambda_{1}(d,\alpha-\gamma,h(t),T)\leq 0$ for $t\geq 0$. Meanwhile, if the assumption $(H)$ holds and the initially occupying habitat is large ($h_0\geq h^{*}$), we also have $\lambda_{1}(d,\alpha-\gamma,h(t),T)\leq 0$. Therefore, it follows
from Lemmas 4.2 and 5.1 that

\noindent\textbf{Theorem 5.1.}
If one of the following assumptions holds:

$(i)$ $F_{h_0}^{+}\neq \emptyset$ and the diffusion is slow ($0<d\leq d_{*}$);

$(ii)$ The assumption $(H)$ holds and the initially occupying habitat is large ($h_0\geq h^{*}$).\\ Then $h_{\infty}=\infty$ and $\lim_{n\rightarrow \infty}u(t+nT,r)=U(t,r)$
locally uniformly for $(t,r)\in[0,T]\times[0,\infty)$, where $U(t,r)$ is the unique positive $T-$periodic solution to $(4.3)$.

The above result implies that, if the average birth rate of a species is greater
than the average death rate, the invasive species with slow diffusion or large
habitat occupation will survive in the new environment. In a biological view,
the species will survive easily in a favorable habitat.

\subsection{Case $\lambda_{1}(d,\alpha-\gamma,h_{0},T)> 0$}

In this subsection, we discuss the case $\lambda_{1}(d,\alpha-\gamma,h_{0},T)> 0$.

\noindent\textbf{Lemma 5.2.}
If $\|u_0(r)\|_{C([0, h_0])}$ is sufficiently large and
one of the following assumptions holds:

$(i)$ The assumption $(H)$ holds and $h_{0}<h^{*}$;

$(ii)$ The diffusion is fast ($d>d^{*}$).\\
Then $h_{\infty}=\infty$.

\noindent\textbf{Proof.}
We fist prove the case $(i)$. Note that
$\lim_{R\rightarrow\infty}\lambda_1(d, \alpha-\gamma, \sqrt{R}, T)<0$ under the assumption $(H)$ from Theorem $3.2$.
Therefore, there exists $R^{*}>0$ such that $\lambda_1(d, \alpha-\gamma, \sqrt{R^{*}}, T)<0$.

Next we construct a suitable lower solution $(\underline{u}(t,r),\underline{h}(t))$ to $(1.1)$, where
$\underline{h}(t)$ satisfies
$\lambda_1(d, \alpha, \gamma,\\ h(R^{*}), T)\leq\lambda_1(d, \alpha-\gamma, \underline{h}(R^{*}), T)\leq\lambda_1(d, \alpha-\gamma, \sqrt{R^{*}}, T)<0$.

Consider the following eigenvalue problem
\begin{align*}
\left\{\begin{array}{l}
\varphi_t-d\Delta\varphi-\frac{1}{2}\varphi_{r}=\lambda\varphi, \quad 0<t<T,~  0<r<1,\\[5pt]
\varphi_{r}(t, 0)=\varphi(t, 1)=0, \quad 0<t<T,\\[5pt]
\varphi(0, r)=\varphi(T, r),\quad 0<r<1.
\end{array}\right.
\tag{5.1}
\end{align*}

It follows from \cite{cc03, hess91} that the above eigenvalue problem has a unique eigenvalue $\lambda_1$, and a corresponding positive eigenfunction $\varphi>0$ in $(t, r)\in [0, T]\times (0, 1)$. By the move-plane argument in \cite{dh94}, we have $\varphi_r(t, r)<0$ in $(t, r)\in [0, T]\times (0, 1]$.

Defining
\begin{align*}
\left\{\begin{array}{l}
\underline{h}(t)=\sqrt{t+\delta}, \quad  t\geq0,\\[5pt]
\underline{u}(t, r)=\frac{M}{(t+\delta)^k}\varphi(\xi, \eta),\quad \xi=\int_{0}^{t}\underline{h}^{-2}(s)ds,\quad \eta=\frac{r}{\sqrt{t+\delta}},\quad t\geq0, \quad 0\leq r\leq\sqrt{t+\delta},
\end{array}\right.
\end{align*}
where $\delta, k, M$ are positive constants to be chosen later, we are now in a position to show that $(\underline{u}, \underline{h})$ is a lower solution of problem $(1.1)$.

It follows from the proof of Theorem $2.1$ that there exists a positive constant $C_{1}$ such that $0\leq u(t,r)\leq C_1$ for $t\in [0, R^{*}], r\in[0, h(t)]$. Under the condition $(1.3)$, there exists a positive constant $L$ such that
$f(t,r,u)=u(\alpha(t,r)-\gamma(t,r)-\beta(t,r)u)\geq-u[\max_{[0,R^{*}]}\alpha_{2}(t)-\min_{[0,R^{*}]}\gamma_{1}(t)-\min_{[0,R^{*}]}\beta_{1}(t)C_{1}]\triangleq -Lu$
for $(t,r)\in [0,R^{*}]\times[0,+\infty)$.

By direct calculations, we obtain
\begin{align*}
&\underline{u}_t-d\Delta\underline{u}-\underline{u}(\alpha(t,r)-\gamma(t,r)-\beta(t,r)\underline{u})\\[5pt]
&=-\frac{M}{(t+\delta)^{k+1}}\{k\varphi(\xi, \eta)
-(t+\delta)[\underline{h}^{-2}(t)\varphi_{\xi}(\xi, \eta)-r\underline{h}^{-2}(t)\underline{h}'(t)\varphi_{\eta}(\xi, \eta)]\\[5pt]
&\quad +d(t+\delta)[\underline{h}^{-2}(t)\varphi_{\eta\eta}(\xi, \eta)+\frac{\underline{h}^{-2}(t)(N-1)}{\eta}\varphi_{\eta}(\xi, \eta)]\\[5pt]
&\quad -(t+\delta)\varphi(\xi, \eta)[\alpha(t,r)-\gamma(t,r)-\beta(t,r)\underline{u}]\}\\[5pt]
&\leq
-\frac{M}{(t+\delta)^{k+1}}\{k\varphi(\xi, \eta)-(t+\delta)\underline{h}^{-2}(t)[\lambda_{1}\varphi(\xi, \eta)+\frac{1}{2}\varphi_{\eta}(\xi, \eta)-r\underline{h}'(t)\varphi_{\eta}(\xi, \eta)]
-L(t+\delta)\varphi(\xi, \eta)    \}\\[5pt]
&=-\frac{M}{(t+\delta)^{k+1}}\{k\varphi(\xi, \eta)-(t+\delta)\underline{h}^{-2}(t)\lambda_1\varphi(\xi, \eta)-L(t+\delta)\varphi(\xi, \eta)\}
\end{align*}
for $0<r<\underline{h}(t), 0<t\leq R^{*}$.

Choosing $0<\delta\leq\min\{1,h_{0}^{2}\}$ and $\lambda_1+L(R^{*}+1)<k$, we obtain
\begin{align*}
\underline{u}_t-d\Delta\underline{u}-\underline{u}(\alpha(t,r)-\gamma(t,r)-\beta(t,r)\underline{u})
\leq-\frac{M}{(t+\delta)^{k+1}}\{k\varphi(\xi, \eta)-\lambda_1\varphi(\xi, \eta)-L(R^{*}+1)\varphi(\xi, \eta)\}<0,
\end{align*}
for $0<r<\underline{h}(t)$ and $0<t\leq R^{*}$.
We may choose may choose $M$ and $\|u_{0}\|_{C([0, h_{0}))}$ being sufficiently large such that
\begin{align*}
\underline{h}'(t)+\mu \underline{u}_{r}(t, h(t))=\frac{1}{2\sqrt{t+\delta}}+\frac{\mu M\varphi_{x}(t, 1)}{(t+\delta)^{k+1/2}}\leq0 \quad \mbox{for} ~0<t\leq R^{*}.
\end{align*}
and
\begin{align*}
\underline{u}(0, r)=\frac{M}{\delta^k}\varphi(0, \frac{r}{\sqrt{\delta}})<u_0(r) \quad\mbox{in} ~[0, \sqrt{\delta}].
\end{align*}
Thus, we have
\begin{align*}
\left\{\begin{array}{l}
\underline{u}_t-d\Delta\underline{u}\leq \underline{u}(\alpha(t,r)-\gamma(t,r)-\beta(t,r)\underline{u}),\quad 0<t\leq R^*, \quad 0<r<\underline{h}(t),\\[5pt]
\underline{u}_{r}(t, 0)=0, \underline{u}(t, \underline{h}(t))=0, \quad 0<t\leq R^{*}, \\[5pt]
\underline{h}'(t)+\mu \underline{u}_{r}(t, \underline{h}(t))\leq0, \quad 0<t\leq R^{*}, \\[5pt]
\underline{u}(0, r)\leq u_0(r), \quad 0\leq r \leq \sqrt{\delta}.
\end{array}\right.
\end{align*}
By the comparison principle to conclude that $\underline{h}(t)\leq h(t)$ in $[0, R^{*}]$. Specially, we derive $h(R^{*})\geq \underline{h}(R^*)=\sqrt{R^*+\delta}\geq\sqrt{R^*}$.
Therefore, $\lambda_1(d, \alpha-\gamma, h(R^{*}), T)\leq \lambda_1(d, \alpha-\gamma, \sqrt{R^{*}}, T)<0$, which together with
Lemma $5.1$ gives that $h_{\infty}=\infty$.

For the case $(ii)$, since $\lim_{d\rightarrow 0}\lambda_{1}(d, \alpha-\gamma, h_{0}, T)<0$ under the assumption $F_{h_{0}}^{+}\neq \emptyset$, then $\lambda_{1}(d, \alpha-\gamma, R, T)<0$ for any $R>h_{0}$ and $0<d<d_{*}$. Thus, we can get the desired result similar to the proof of $(i)$.\quad$\Box$

\noindent\textbf{Lemma 5.3.} $h_{\infty}<\infty$
if $\|u_{0}\|_{C([0,h_{0}])}$ is sufficiently small.

\noindent\textbf{Proof.}
Recall that $\lambda_1(d, \alpha-\gamma, h_0, T)$ is the principal eigenvalue
of the following problem
\begin{align*}
\left\{\begin{array}{l}
\varphi_t-d\Delta\varphi=\varphi(\alpha(t,r)-\gamma(t,r))+\lambda\varphi, \quad \mbox{in}~[0,T]\times B_{h_{0}},\\[5pt]
\varphi(t, r)=0, \quad \mbox{on}~[0,T]\times \partial B_{h_{0}},\\[5pt]
\varphi(0, r)=\varphi(T, r),\quad \mbox{in}~B_{h_{0}},
\end{array}\right.
\tag{5.2}
\end{align*}
and $\varphi\in C^{1,2}([0,T]\times\bar{B}_{h_{0}})$ with
$\|\varphi\|_{L^{\infty}([0,T]\times B_{h_{0}})}=1$ is the corresponding positive $T-$periodic eigenfunction.

We first prove a fact that there exists some $M_2>0$ such that $\varphi_{r}(t, r)\leq M_2\varphi(t, r)$,
$\forall (t, r)\in [0, T]\times [0, h_0]$.

Let $\psi(t, r)=\varphi_{r}(t, r)-M_{1}(2h_{0}-r)\varphi(t, r)$, where $M_{1}$ is a positive constant which will be
determined later. Direct calculation gives that for $\forall (t, r)\in [0, T]\times [0, h_0]$,
\begin{align*}
&\psi_{t}-d(\psi_{rr}+\frac{N-1}{r}\psi_r)=\varphi_{tr}-M_1(2h_0-r)\varphi_t
-d\{\varphi_{rrr}-M_1(2h_0-r)\varphi_{rr}+2M_1\varphi_{r}\\[5pt]
&\quad+\frac{N-1}{r}[\varphi_{rr}-M_1(2h_0-r)\varphi_{r}+M_1\varphi]\}\\[5pt]
&=\varphi_{tr}-M_1(2h_0-r)\varphi_t-d(\varphi_{rrr}+\frac{N-1}{r}\varphi_{rr})
-d\{-M_1(2h_0-r)\varphi_{rr}+2M_1\varphi_{r}\\[5pt]
&\quad+\frac{N-1}{r}[-M_1(2h_0-r)\varphi_{r}+M_1\varphi]\}\\[5pt]
&=-M_1(2h_0-r)\varphi_t+(\alpha(t, r)-\gamma(t, r))_{r}\varphi+(\alpha(t, r)-\gamma(t, r)+\lambda_1)\varphi_{r}
-\frac{N-1}{r^2}d\varphi_{r}\\[5pt]
&\quad+dM_1(2h_0-r)(\varphi_{rr}+\frac{N-1}{r}\varphi_{r})-2dM_1[\psi+M_1(2h_0-r)\varphi]
-\frac{N-1}{r}dM_1\varphi\\[5pt]
&=(\alpha(t, r)-\gamma(t, r))_{r}\varphi+(\alpha(t, r)-\gamma(t, r)+\lambda_1)\varphi_{r}-\frac{N-1}{r^2}d\varphi_{r}\\[5pt]
&\quad-M_1(2h_0-r)(\alpha(t, r)-\gamma(t, r)+\lambda_1)\varphi
-2dM_1[\psi+M_1(2h_0-r)\varphi]-\frac{N-1}{r}dM_1\varphi\\[5pt]
&=[(\alpha(t, r)-\gamma(t, r))_{r}-2dM_1^2(2h_0-r)]\varphi+(\alpha(t, r)-\gamma(t, r)+\lambda_1)\psi-\frac{N-1}{r^2}d\varphi_{r}\\[5pt]
&\quad-2dM_1\psi-\frac{N-1}{r}dM_1\varphi\\[5pt]
&=[(\alpha(t, r)-\gamma(t, r))_{r}-2dM_1^2(2h_0-r)]\varphi+(\alpha(t, r)-\gamma(t, r)+\lambda_1-2dM_1)\psi-\frac{N-1}{r}dM_1\varphi\\[5pt]
&\quad-\frac{N-1}{r^2}d[\psi+M_1(2h_0-r)\varphi]\\[5pt]
&=[(\alpha(t, r)-\gamma(t, r))_{r}-2dM_1^2(2h_0-r)]\varphi+(\alpha(t, r)-\gamma(t, r)+\lambda_1-2dM_1-\frac{N-1}{r^2}d)\psi\\[5pt]
&\quad-\frac{N-1}{r}dM_1\varphi-\frac{N-1}{r^2}dM_1(2h_0-r)\varphi\\[5pt]
&\leq [(\alpha(t, r)-\gamma(t, r))_{r}-2dM_1^2(2h_0-r)]\varphi+(\alpha(t, r)-\gamma(t, r)+\lambda_1-2dM_1-\frac{N-1}{r^2}d)\psi\\[5pt]
&\leq (\alpha(t, r)-\gamma(t, r)+\lambda_1-2dM_1-\frac{N-1}{r^2}d)\psi,
\end{align*}
provided that $M_1\geq\sqrt{\frac{\|\alpha-\gamma\|_{C^{0, 1}([0, T]\times[0, h_0])}}{2dh_0}}$. Clearly,
$\psi(t, 0)=\varphi_{r}'(t, 0)-2M_1h_0\varphi(t, 0)\leq 0$ and $\psi(t, h_0)=\varphi_{r}(t, h_0)-M_1h_0\varphi(t, h_0)
=\varphi_{r}(t, h_0)<0$ for $0\leq t \leq T$. Moreover, by the Hopf boundary lemma for parabolic equations, we have $\varphi_{r}(0, h_0)<0$.
Thus if we take $M_1\geq\frac{\|\varphi\|_{C^{1, 2}([0, T]\times\bar{B}_{h_0})}}{2h_0}$, we can derive that
$\psi(0, r)=\varphi_{r}(0, r)-M_1(2h_0-r)\varphi(0, r)\leq 0$ for $0\leq r \leq h_0$. If we further choose
$$
M_1=\max\left\{\frac{\|\alpha-\gamma\|_{C([0, T]\times[0, h_0])+\lambda_1}}{2d}, \sqrt{\frac{\|\alpha-\gamma\|_{C^{0, 1}([0, T]\times[0, h_0])}}{2dh_0}}, \frac{\|\varphi\|_{C^{1, 2}([0, T]\times[0, h_0])}}{2h_0}\right\},
$$
then by the maximum principle we see that $\psi(t, r)\leq 0$ for all $(t, r)\in[0, T]\times[0, h_0]$, i.e.,
$\varphi_{r}(t, r)\leq M_2\varphi(t, r)$ for all $(t, r)\in[0, T]\times[0, h_0]$ with $M_2=h_0M_1$.

Next, we will construct a suitable supersolution to $(1.1)$. For $t>0$ and $r\in [0,\bar{h}(t)]$, we define
\begin{align*}
\bar{h}(t)=h_{0}\tau(t),\quad \tau(t)=1+\varepsilon-\frac{\varepsilon}{2}e^{-\varepsilon t}, \\[5pt]
\bar{u}(t,r)=Me^{-\varepsilon t}\varphi(\int_{0}^{t}\tau^{-2}(s)ds, \frac{h_{0}}{\bar{h}(t)}r),
\end{align*}
where $\varepsilon, M$ are positive constants to be chosen later, and $\bar{u}(t,r)$ is well defined for all $t>0$.
For simplicity, denote $\xi=\int_{0}^{t}\tau^{-2}(s)ds$, $\eta=r\tau^{-1}(t)$. Thus, $\bar{u}(t,r)=Me^{-\varepsilon t}\varphi(\xi,\eta)$.
By direct calculations yield
\begin{align*}
&\bar{u}_t-d\Delta\bar{u}-\bar{u}(\alpha(t,r)-\gamma(t,r)-\beta(t,r)\bar{u})\\[5pt]
&=Me^{-\varepsilon t}\{-\varepsilon\varphi(\xi, \eta)
-r\tau^{-2}(t)\tau'(t)\varphi_{\eta}(\xi, \eta)+\tau^{-2}(t)\varphi_{\xi}(\xi, \eta)
-d\tau^{-2}(t)\varphi_{\eta\eta}(\xi, \eta)\\[5pt]
&\quad -\frac{d\tau^{-2}(t)(N-1)}{\eta}\varphi_{\eta}(\xi, \eta)
-\varphi(\xi, \eta)[\alpha(t,r)-\gamma(t,r)-\beta(t,r)Me^{-\varepsilon t}\varphi(\xi, \eta)]\}\\[5pt]
&=Me^{-\varepsilon t}\{-\varepsilon \varphi(\xi, \eta)
+\tau^{-2}(t)[\lambda_{1}\varphi(\xi,\eta)+\varphi(\xi,\eta)(\alpha(\xi,\eta)-\gamma(\xi,\eta))]
\\[5pt]
&\quad-\gamma\varphi_{\eta}\tau^{-2}(t)\tau^{\prime}(t)
 -\varphi(\xi, \eta)[\alpha(t,r)-\gamma(t,r)-\beta(t,r)Me^{-\varepsilon t}\varphi(\xi, \eta)]\}  \\[5pt]
&\geq Me^{-\varepsilon t}\varphi(\xi, \eta)
\{-\varepsilon+\tau^{-2}(t)[\lambda_{1}+\alpha(\xi, \eta)-\gamma(\xi, \eta)]
 -[\alpha(t,r)-\gamma(t,r)]-h_0\varepsilon^2M_2\}
\\[5pt]
&=\bar{u}(\xi,\eta)\{-\varepsilon+\tau^{-2}(t)\lambda_{1}+
\tau^{-2}(t)[\alpha(\xi, \eta)-\gamma(\xi, \eta)]-[\alpha(t,r)-\gamma(t,r)]-h_0\varepsilon^2M_2\}.
\end{align*}
It easy to observe that for $t>0$
\begin{align*}
1+\frac{\varepsilon}{2}\leq \tau(t)\leq 1+\varepsilon,\quad
h_{0}(1+\frac{\varepsilon}{2})\leq \bar{h}(t)\leq h_{0}(1+\varepsilon).
\end{align*}
Therefore,
\begin{align*}
(1+\varepsilon)^{-2}t\leq \xi\leq(1+\frac{\varepsilon}{2})^{-2}t,\quad
(1+\varepsilon)^{-1}r\leq \eta\leq(1+\frac{\varepsilon}{2})^{-1}r.
\end{align*}
Moreover, one can calculate
\begin{align*}
&\left|\tau^{-2}(t)[\alpha(\xi,\eta)-\gamma(\xi,\eta)]-[\alpha(t,r)-\gamma(t,r)]\right|\\[5pt]
&\leq \tau^{-2}(t)\left|[\alpha(\xi,\eta)-\gamma(\xi,\eta)]-[\alpha(t,r)-\gamma(t,r)]\right|
 +\left|(\tau^{-2}(t)-1)[\alpha(t,r)-\gamma(t,r)]\right|\\[5pt]
&\leq \left|[\alpha(\xi,\eta)-\gamma(\xi,\eta)]-[\alpha(t,r)-\gamma(t,r)]\right|
 +\left|(\tau^{-2}(t)-1)\right|\|\alpha(t,r)-\gamma(t,r)\|_{C([0,T]\times [0,2h_{0}])}\\[5pt]
&\leq \|\alpha-\gamma\|_{C^{\frac{\nu_{0}}{2},\nu_{0}}([0,T]\times [0,2h_{0}])}
 \left[|\xi-t|]^{\frac{\nu_{0}}{2}}+|\eta-r|^{\nu_{0}}\right]
 +\left|(\tau^{-2}(t)-1)\right|\|\alpha(t,r)-\gamma(t,r)\|_{C([0,T]\times [0,2h_{0}])}.
\end{align*}
Since $\tau^{-2}(t)\rightarrow1$, $\xi\rightarrow t$, $\eta\rightarrow r$ as $\varepsilon\rightarrow 0$,
the above inequality implies that there exists $\varepsilon_{1}>0$ small such that
\begin{align*}
\left|\tau^{-2}(t)[\alpha(\xi,\eta)-\gamma(\xi,\eta)]-[\alpha(t,r)-\gamma(t,r)]\right|
\leq \frac{1}{4}\lambda_{1},
\quad \mbox{for}~\varepsilon\leq \varepsilon_{1}.
\end{align*}
Obviously, we can choose some $\varepsilon_{2}>0$ small such that
\begin{align*}
\frac{1}{(1+\varepsilon)^{2}}\geq \frac{3}{4}, \quad \mbox{and}~h_0\varepsilon^2M_2\leq\frac{1}{4}\lambda_1, \quad \mbox{for}~\varepsilon\leq \varepsilon_{2}.
\end{align*}
Set $\varepsilon=\min\{1,\frac{1}{4}\lambda_{1},\varepsilon_{1},\varepsilon_{2}\}>0$,
then for $t>0$ and $0\leq x\leq \bar{h}(t)$, we get
\begin{align*}
\bar{u}_t-d\Delta\bar{u}-\bar{u}(\alpha(t,r)-\gamma(t,r)-\beta(t,r)\bar{u})
\geq \bar{u}(\xi,\eta)\left\{-\frac{1}{4}\lambda_{1}+\frac{3}{4}\lambda_{1}
-\frac{1}{4}\lambda_{1}-\frac{1}{4}\lambda_{1}\right\}
=0.
\end{align*}
Moreover, we now choose $M>0$ sufficiently large and $\|u_{0}\|_{C([0,h_{0}])}$ is sufficiently small such that
\begin{align*}
u_{0}(r)\leq M\varphi(0, \frac{r}{1+\frac{\varepsilon}{2}})=\bar{u}(0,r)
\quad \mbox{for}~r\in [0,h_{0}].
\end{align*}
Then
\begin{align*}
&\bar{h}(t)=\frac{1}{2}h_{0}\varepsilon^{2}e^{-\varepsilon t},\\[5pt]
&-\mu\bar{u}(t,\bar{h}(t))=\mu M e^{-\varepsilon t}
 \left|\varphi_{\eta}(\int_{0}^{t}\tau^{-2}(s)ds,h_{0})\right|\tau^{-1}(t)
\leq \frac{c_{0}\mu M}{1+\frac{\varepsilon}{2}}e^{-\varepsilon t},
\end{align*}
where $c_{0}=\max_{t\in[0,\infty)}\left|\varphi_{\eta}(\int_{0}^{t}\tau^{-2}(s)ds,h_{0})\right|$.
Thus, if we choose
\begin{align*}
\mu_{0}\leq \frac{\varepsilon^{2}(1+\frac{\varepsilon}{2})h_{0}}{2Mc_{0}},
\end{align*}
then $\bar{h}(t)\geq-\mu\bar{u}_{r}(t,\bar{h}(t))$ for $0<\mu\leq\mu_{0}$,
and thus $(\bar{u}(t,r),\bar{h}(t))$ satisfies
\begin{align*}
\left\{\begin{array}{l}
\bar{u}_t-d\Delta\bar{u}\geq \bar{u}(\alpha(t,r)-\gamma(t,r)-\beta(t,r)\bar{u}),\quad t>0, 0<r<\bar{h}(t),\\[5pt]
\bar{u}_{r}(t, 0)=0, \bar{u}(t, \bar{h}(t))=0, \quad t>0, \\[5pt]
\bar{h}^{\prime}(t)+\mu \bar{u}_{r}(t, \bar{h})\geq0, \quad t>0, \\[5pt]
\bar{u}(0, r)\geq u_0(r), \quad 0\leq r \leq \bar{h}(0).
\end{array}\right.
\end{align*}
By the comparison principle, we obtain that $h(t)\leq \bar{h}(t)$ and
$u(t,r)\leq \bar{u}(t,r)$ for $0\leq r\leq h(t)$ and $t>0$. This implies that
$h_{\infty}\leq \lim_{t\rightarrow \infty}\bar{h}(t)=h_{0}(1+\varepsilon)<\infty$.\quad$\Box$

In fact, by using the same arguments as above, we can prove a more general result as follows.

\noindent\textbf{Lemma 5.4.} For any given $u_{0}$ satisfying $(1.2)$, there exists $\mu_{0}>0$
depending on $u_{0}$ such that $h_{\infty}<\infty$ provided $0<\mu\leq\mu_{0}$.

In the following, with the parameters $h_{0}$ and $d$ satisfying $0<h_{0}<h^{*}$ and $d>d^{*}$, respectively, which implies that $\lambda_{1}(d, \alpha-\gamma, h_{0}, T)>0$, and $u_{0}$ fixed, let us discuss the effect of the coefficient $\mu$ on the spreading and vanishing when $\mu$ is sufficiently large. 

\noindent\textbf{Lemma 5.5.} There exists $\mu^{0}>0$ depending on $u_{0}, h_{0}$ and $d$ such that $h_{\infty}=\infty$
if $\mu>\mu^{0}$.

\noindent\textbf{Proof.} The idea of this proof comes from Lemma 3.6 in \cite{pz13}.
Note that from $(4.4)$ there  exists a positive constant $M$ such that $0<u(t,r)\leq M$
for $t\in (0,\infty)$ and $0\leq r<h(t)$, where $u(t,r)$ is the unique positive solution to
problem $(1.1)$. Thus, there exists a positive constant $\bar{L}$ such that
\begin{align*}
u(\alpha(t,r)-\gamma(t,r)-\beta(t,r)u)
\geq u[-\max_{[0,T]}\gamma_{2}(t)-\max_{[0,T]}\beta_{2}(t)M]\triangleq -\bar{L}u,
\end{align*}
for all $t\in(0,\infty)$ and $r\in [0,h(t)]$. We next consider the auxiliary free boundary
problem
\begin{align*}
\left\{\begin{array}{l}
W_t-d\Delta W=-\bar{L}W,\quad t>0, \quad 0<r<\delta(t),\\[5pt]
W_{r}(t, 0)=0, W(t, \delta(t))=0, \quad t>0, \\[5pt]
\delta^{\prime}(t)+\mu W_{r}(t, \delta(t))=0, \quad t>0, \\[5pt]
W(0, r)= u_0(r), \quad 0\leq r \leq \delta(0)=h_{0}.
\end{array}\right.
\tag{5.3}
\end{align*}
Arguing as in proving the existence and uniqueness of the solution to problem $(1.1)$, one will easily see that $(5.3)$
also admits a unique global solution $(W,\delta)$ and $\delta^{\prime}(t)>0$ for $t>0$. To stress the dependence of the solutions
on the parameter $\mu$, in the sequel, we always write $(u^{\mu}, h^{\mu})$ and $(W^{\mu}, \delta^{\mu})$ instead of $(u,h)$
and $(W, \delta)$. By the comparison principle, we have
\begin{align*}
u^{\mu}(t,r)\geq W^{\mu}(t,r), \quad h^{\mu}(t)\geq \delta^{\mu}(t)\quad \mbox{for any}~t\geq 0, r\in [0, \delta^{\mu}(t)].
\tag{5.4}
\end{align*}
In what follows, we are in a position to prove that for all large $\mu$,
\begin{align*}
\delta^{\mu}(T^{*})\geq T^{*}h^{*},
\tag{5.5}
\end{align*}
for any positive constant $T^{*}>1$.

To this end, we first choose a smooth function $\underline{\delta}(t)$ with $\underline{\delta}(0)=\frac{h_{0}}{T^{*}}$,
$\underline{\delta}^{\prime}(t)>0$ and $\underline{\delta}(T^{*})=T^{*}h^{*}$. We then consider the following initial-boundary
problem
\begin{align*}
\left\{\begin{array}{l}
\underline{W}_t-d\Delta \underline{W}=-\bar{L}\underline{W},\quad t>0, \quad 0<r<\underline{\delta}(t),\\[5pt]
\underline{W}_{r}(t, 0)=0, \underline{W}(t, \underline{\delta}(t))=0, \quad t>0, \\[5pt]
\underline{W}(0, r)= \underline{W}_0(r), \quad 0\leq r \leq \frac{h_{0}}{T^{*}}.
\end{array}\right.
\tag{5.6}
\end{align*}
Here, for the smooth initial value $\underline{W}_0(r)$, we require
\begin{align*}
0<\underline{W}_0(r)\leq u_{0}(r)\quad \mbox{on}~[0,\frac{h_{0}}{T^{*}}], \quad
\underline{W}_0^{\prime}(0)=\underline{W}_0(\frac{h_{0}}{T^{*}})=0, \quad\underline{W}_0^{\prime}(\frac{h_{0}}{T^{*}})<0. \tag{5.7}
\end{align*}
The standard theory for parabolic equations ensures that $(5.6)$ has a unique positive solution $\underline{W}$,
and $\underline{W}_{r}(t,\underline{\delta}(t))<0$ for all $t\in [0,T^{*}]$ due to the well-known Hopf boundary lemma.
According to our choice of $\underline{\delta}(t)$ and $\underline{W}_0(r)$, there is a constant $\mu^{0}>0$ such that
for all $\mu\geq \mu^{0}$,
\begin{align*}
\underline{\delta}^{\prime}(t)\leq -\mu\underline{W}_{r}(t,\underline{\delta}(t)), \quad \mbox{for all}~0\leq t\leq T^{*}. \tag{5.8}
\end{align*}

On the other hand, for problem $(5.3)$, we can establish the comparison principle analogous with lower solution to Lemma 2.1 by the same argument. Thus, note that $\underline{\delta}(0)=\frac{h_{0}}{T^{*}}<\delta^{\mu}(0)$, it follows from
$(5.3)$ and $(5.6)-(5.8)$ that
\begin{align*}
W^{\mu}(t,r)\geq \underline{W}(t,r), \quad \delta^{\mu}(t)\geq \underline{\delta}(t),
\quad t\in[0,T^{*}], r\in[0,\underline{\delta}(t)],
\end{align*}
which implies that $\delta^{\mu}(T^{*})\geq \underline{\delta}(T^{*})=T^{*}h^{*}$,
and so $(5.5)$ holds. Thus, in view of $(5.4)$ and $(5.5)$,
we find $h_{\infty}=\lim_{t\rightarrow+\infty}h^{\mu}(t)>h^{\mu}(T^{*})\geq T^{*}h^{*}>h^{*}$.
This, together with Lemma 4.1, yields the desired result. \quad$\Box$

Using the parameter $\mu$, we will give the sharp criteria for spreading and vanishing for $\lambda_{1}(d,\alpha,\gamma,\\h_{0},T)>0$.

\noindent\textbf{Theorem 5.2.}
There exists $\mu^{*}>0$, depending on $u_0$, $h_0$ and $d$, such that
$h_{\infty}=\infty$ if $\mu>\mu^{*}$, and $h_{\infty}\leq h^{*}$ if $\mu\leq \mu^{*}$.

\noindent\textbf{Proof.}
Recall Lemma 5.4 and 5.5, the proof is similar to that of Theorem 2.10 in \cite{dg11}, so we omit the details.\quad$\Box$

\subsection{Sharp criteria}

In this subsection, we will establish the sharp criteria by selecting $d, h_0, \mu$ and $u_0$ as varying factors to distinguish the spreading and vanishing for the invasive species.

If $h_{0}$ is fixed, the spreading or vanishing of an invasive species
depends on the diffusion rate $d$, the initial number $u_{0}(r)$ of the species and the expanding ability $\mu$.
Now we give the sharp criteria for species spreading and vanishing.
The result directly follows from Theorems $5.1-5.2$, and Lemmas $5.2-5.3$.\\

\noindent\textbf{Theorem 5.3.} Suppose $(H)$ holds. \\
$(i)$ If $F_{h_{0}}^{+}\neq\emptyset$, then there exists $d_{*}\in(0,\infty)$
such that the spreading happens for any initial value if the diffusion is slow
($0<d\leq d_{*}$);\\
$(ii)$ There exists $d^{*}\in(0,\infty)$ such that when the diffusion is fast ($d>d^{*}$), we have

$(ii.1)$ one of the following results holds: vanishing occurs if $\|u_{0}\|_{C([0,h_{0}])}$ is sufficiently small;
spreading occurs if $\|u_{0}\|_{C([0,h_{0}])}$ is sufficiently large.

$(ii.2)$ there exists $\mu^{*}>0$, depending on $u_0$, $h_0$ and $d$, such that spreading happens if
$\mu>\mu^{*}$, and vanishing occurs if $\mu\leq\mu^{*}$;\\
$(iii)$ moreover, $\mu^*=0$ if the diffusion is slow ($0<d\leq d_{*}$).

Theorem 5.3 implies that the invasive species with slow diffusion will spread
in the total habitat, while if the ability of migratory is big, the vanishing or spreading
of the species in the new environment is determined by the initial number and the expanding
capability.

Similarly, if $d$ is fixed, the spreading or vanishing of an invasive species
depends on the radius $h_{0}$ of the initial occupying habitat $B_{h_{0}}$, the initial number $u_{0}(r)$ of
the species and the expanding capability $\mu$. The following result obviously derives
from Theorems $5.1-5.2$, and Lemmas $5.2-5.3$.

\noindent\textbf{Theorem 5.4.} Suppose $(H)$ holds. There exist $h^{*}\in(0,\infty)$ such that\\
$(i)$ the spreading happens for any initial value if $h_{0}\geq h^{*}$;\\
$(ii)$ when $h_{0}<h^{*}$, then

$(ii.1)$  one of the following results holds: vanishing occurs if $\|u_{0}\|_{C([0,h_{0}])}$ is sufficiently small;
spreading occurs if $\|u_{0}\|_{C([0,h_{0}])}$ is sufficiently large.

$(ii.2)$ there exists $\mu^{*}>0$, depending on $u_0$, $h_0$ and $d$, such that spreading happens if
$\mu>\mu^{*}$, and vanishing occurs if $\mu\leq\mu^{*}$;\\
$(iii)$ moreover, $\mu^*=0$ if $h_0\geq h^*$.

Theorem 5.4 tells us an unfavorable habitat is bad for the species with small number at
the beginning, the endangered species in an unfavorable habitat will become extinct in the
future. However, even the habitat is unfavorable, if the occupying area $B_{h_{0}}$ is
beyond a critical size, namely $h_{0}\geq h^{*}$, then regardless of the initial population
size $u_{0}(r)$, spreading always happens (see $(i)$ in Theorem 5.4). And if $h_{0}<h^{*}$,
spreading is also possible for big initial population size $u_{0}(r)$.

If $d$ and $h_0$ are fixed, the initial number $u_0(r)$ governs the spreading and vanishing of the invasive species.

\noindent\textbf{Theorem 5.5.}
Assume that $(H)$ holds. Let $(u(t, r), h(t))$ be the solution of problem $(1.1)$ with $u_0(r)=\sigma\zeta(r)$ for some positive constant $\sigma$. Then there exists $\sigma_0\triangleq\sigma_0(d, h_0)\in[0, +\infty)$ such that spreading occurs for $\sigma>\sigma_0$ and vanishing happens for $0<\sigma\leq\sigma_0$.

\noindent\textbf{Proof.}
The result follows from Theorem 3.2, Lemmas $5.1-5.3$. The proof is similar to Theorem 5.7 in \cite{llz14} with obvious modification. So we omit it.

Theorems $5.1-5.5$ imply that slow diffusion, large occupying habitat and big initial
population number are benefit for the species to survive in the new environment.

\section{Spreading speed}

In this section, we always suppose $(H)$ holds. In the spreading case, we will give the asymptotic spreading speed of the free boundary $x=h(t)$.

Consider the following problem
\begin{align*}
\left\{\begin{array}{l}
U_t-d\Delta U+k(t)U_r=U(a(t)-b(t)U), \quad (t, r)\in[0, T]\times(0, \infty),\\[5pt]
U(t, 0)=0,\quad t\in[0, T],\\[5pt]
U(0, r)=U(T, r), \quad r\in (0, \infty),
\end{array}\right.
 \tag{6.1}
\end{align*}
where $d>0$ is a given constant, and $k, a, b$ are given T-periodic H\"{o}lder continuous functions with $a, b$ positive and $k$ nonnegative.

It follows from Proposition 2.1, 2.3 and Theorem 2.4 in \cite{dgp13} that

\noindent\textbf{Lemma 6.1.}
For any given positive T-periodic functions $a, b\in C^{\frac{\nu_0}{2}}([0, T])$ and any nonnegative continuous T-periodic function $k\in C^{\frac{\nu_0}{2}}([0, T])$, problem $(6.1)$ admits a positive T-periodic solution $U^k\in C^{1, 2}([0, T]\times[0, \infty))$ if and only if $\frac{1}{T}\int_{0}^{T}a(t)dt>\frac{1}{T^2}(\int_{0}^{T}k(t)dt)^2/(4d)$. Moreover, either $U^k\equiv0$ or $U^k>0$ in $[0, T]\times[0, \infty)$. Furthermore, if $U^k>0$, then it is the only positive solution of problem $(6.1)$, $U^k_{r}(t, r)>0$ in $[0, T]\times[0, \infty)$ and $U^k(t, r)\rightarrow V(t)$ uniformly for $t\in[0, T]$ as $r\rightarrow+\infty$, where $V(t)$ is the unique positive solution of the problem
\begin{align*}
\left\{\begin{array}{l}
\frac{dV}{dt}=V(a(t)-b(t)V), \quad t\in[0, T],\\[5pt]
V(0)=V(T).
\end{array}\right.
 \tag{6.2}
\end{align*}
In addition, for any given nonnegative T-periodic function $k_1\in C^{\frac{\nu_0}{2}}([0, T])$, the assumption $k_1\leq,\not\equiv k$ implies $U^{k_1}_{r}(t, 0)>U^k_{r}(t, 0), U^{k_1}(t, r)>U^k(t, r)$ for $(t, r)\in[0, T]\times(0, +\infty)$. Besides, for each $\mu>0$, there exists a positive continuous T-periodic function $k_0(t)=k_0(\mu, a, b)(t)>0$ such that $\mu U^{k_0}_{r}(t, 0)=k_0(t)$ on $[0, T]$. Further, $0<\frac{1}{T}\int_{0}^{T}k_0(\mu, a, b)(t)dt<2\sqrt{\frac{d}{T}\int_{0}^{T}a(t)dt}$ for every $\mu>0$.

By the assumption $(H)$ and the condition $(1.3)$, we have that, for any $\varepsilon>0$, there exists a $R_{*}\triangleq R(\varepsilon)>1$ such that for $r\geq R_{*}$,
\begin{align*}
&\alpha(t, r)-\gamma(t, r)\leq\eta^{*, \varepsilon}\triangleq\eta^*(t)+\varepsilon,\\[5pt]
&\alpha(t, r)-\gamma(t, r)\geq\eta_{*, \varepsilon}(t)\triangleq\eta_{*}(t)-\varepsilon,\\[5pt]
&\beta(t, r)\leq\beta_2^{\varepsilon}(t)\triangleq\beta_2(t)+\varepsilon,\\[5pt]
&\beta(t, r)\geq\beta_{1, \varepsilon}(t)\triangleq\beta_1(t)-\varepsilon,
\end{align*}
and $\eta_{*, \varepsilon}(t), \eta^{*, \varepsilon}, \beta_{1, \varepsilon}(t), \beta_2^{\varepsilon}(t)$ are T-periodic H\"{o}lder continuous functions.

We are now ready to give our result.

\noindent\textbf{Theorem 6.1.}
If $h_{\infty}=\infty$, then
\begin{align*}
\limsup_{t\rightarrow+\infty}\frac{h(t)}{t}\leq\frac{1}{T}\int_{0}^{T}k_0(\mu, \eta^*, \beta_1)(t)dt,
 \tag{6.3}
\end{align*}
\begin{align*}
\liminf_{t\rightarrow+\infty}\frac{h(t)}{t}\geq\frac{1}{T}\int_{0}^{T}k_0(\mu, \eta_*, \beta_2)(t)dt,
 \tag{6.4}
\end{align*}
where $k_0(\mu, \cdot, \cdot)$ is given in Lemma 6.1.

\noindent\textbf{Proof.}
The proof is a simple modification of that for Theorem 4.4 in \cite{dgp13}. So we only briefly describe the main steps.

Step 1. The unique positive solution $\hat{U}$ of $(4.3)$ satisfies
\begin{align*}
\limsup_{r\rightarrow+\infty}\hat{U}(t, r)\leq\bar{V}(t),\quad
\liminf_{r\rightarrow+\infty}\hat{U}(t, r)\geq\underline{V}(t),\quad \mbox{for}~t\in[0, T],
 \tag{6.5}
\end{align*}
where $\underline{V}(t)$ and $\bar{V}(t)$ are the unique positive T-periodic solutions of
\begin{align*}
\left\{\begin{array}{l}
\frac{d\bar{V}}{dt}=\bar{V}(\eta^*(t)-\beta_1(t)\bar{V}), \quad t\in[0, T],\\[5pt]
\bar{V}(0)=\bar{V}(T),
\end{array}\right.
 \tag{6.6}
\end{align*}
and
\begin{align*}
\left\{\begin{array}{l}
\frac{d\underline{V}}{dt}=\underline{V}(\eta_*(t)-\beta_2(t)\underline{V}), \quad t\in[0, T],\\[5pt]
\underline{V}(0)=\underline{V}(T),
\end{array}\right.
 \tag{6.7}
\end{align*}
respectively.

This step can be obtained by the sub-supersolution argument, which is similar to the step 1 in Theorem 4.4 in \cite{dgp13} with the help of Proposition 4.1 and 4.2.

Step 2. We construct two suitable supersolution and subsolution to problem $(1.1)$ to prove $(6. 3)$ and $(6.4)$, respecyively.

By $(6.5)$, there exists $\hat{R}^*\triangleq \hat{R}^*(\varepsilon)>R_{*}>1$ such that
\begin{align*}
\underline{V}_{\frac{\varepsilon}{2}}(t)\leq \hat{U}(t, r)\leq\bar{V}_{\frac{\varepsilon}{2}}(t),
\quad \forall (t, r)\in[0, T]\times[\hat{R}^*, +\infty),
\end{align*}
where $\underline{V}_{\frac{\varepsilon}{2}}(t)$ and $\bar{V}_{\frac{\varepsilon}{2}}(t)$
are the unique positive T-periodic solution to $(6.6)$ and $(6.7)$ by substituting $\eta^*, \beta_1, \eta_*, \beta_2$ with $\eta^{*, \frac{\varepsilon}{2}}, \beta_{1, \frac{\varepsilon}{2}}, \eta_{*, \frac{\varepsilon}{2}}, \beta_2^{\frac{\varepsilon}{2}}$, respectively.

Since $h_{\infty}=+\infty$ and $\lim_{n\rightarrow+\infty}u(t+nT, r)=\hat{U}(t, r)$, there exists a positive integer $K=K(\hat{R}^*)$ such that with $\hat{T}\triangleq KT, h(\hat{T})>3\hat{R}^*$ and $u(t+\hat{T}, 2\hat{R}^*)<\bar{V}_{\varepsilon}(t)$ for all $t>0$.

Let $\hat{u}(t, r)=u(t+\hat{T}, r+2\hat{R}^*)$ and $\hat{h}(t)=h(t+\hat{T})-2\hat{R}^*$.
Direct calculations yield
\begin{align*}
\left\{\begin{array}{l}
\hat{u}_t-d(\hat{u}_{rr}+\frac{N-1}{r+2\hat{R}^*}\hat{u}_{r})=u(\alpha(t+\hat{T}, r+2\hat{R}^*)-\gamma(t+\hat{T}, r+2\hat{R}^*)\\[5pt]
\qquad\qquad\qquad\qquad\qquad\quad-\beta(t+\hat{T}, r+2\hat{R}^*)\hat{u}),
t>0,\quad 0<r<\hat{h}(t),\\[5pt]
\hat{u}(t, 0)=u(t+\hat{T}, 2\hat{R}^*), \hat{u}(t, \hat{h}(t))=0,\quad t>0,\\[5pt]
\hat{h}'(t)=-\mu \hat{u}_{r}(t, \hat{h}(t)), \quad t>0,\\[5pt]
\hat{u}(0, r)=u(\hat{T}, r+2\hat{R}^*), \quad 0\leq r\leq \hat{h}_{0}.
\end{array}\right.
 \tag{6.8}
\end{align*}
It follows from the comparison principle that
\begin{align*}
\hat{u}(t, r)\leq w(t)\leq \bar{V}_{\varepsilon}(t)(1-\varepsilon)^{-1},\quad \mbox{for}~t\geq\hat{T}\quad \mbox{and}~0<r<\hat{h}(t),
\end{align*}
where $w(t)$ is the unique positive solution of the problem
\begin{align*}
\left\{\begin{array}{l}
\frac{dw}{dt}=w(\eta^{*, \varepsilon}(t)-\beta_{1, \varepsilon}(t)w), \quad t>0,\\[5pt]
w(0)=\max\left\{\bar{V}_{\varepsilon}, \|\hat{u}(0, \cdot)\|_{\infty}\right\}.
\end{array}\right.
\end{align*}

On the other hand, replacing $\alpha(t+\hat{T}, r+2\hat{R}^*)-\gamma(t+\hat{T}, r+2\hat{R}^*)$ and $\beta(t+\hat{T}, r+2\hat{R}^*)$ by $\eta_{*, \varepsilon}(t)$ and $\beta_2^{\varepsilon}$, respectively. We can construct a lower solution to $(6.8)$, then we obtain that $\hat{u}$ has a subbound. That is
\begin{align*}
\liminf_{r\rightarrow+\infty}\hat{u}(t+nT, r)\geq \underline{V}_{\varepsilon}(t)\quad\mbox{locally uniformly for}~(t, r)\in[0, T]\times[0, \infty).
\end{align*}

Now we use the sub-supersolution argument to get our desired result.
Let $U_{\varepsilon}=U_{\eta^{*, \varepsilon}, \beta_{1, \varepsilon}, k^{\varepsilon}}$ and
$V_{\varepsilon}=U_{\eta^{*, \varepsilon}, \beta_2^{\varepsilon}, k_{\varepsilon}}$
be the unique positive solutions to $(6.1)$ with $a(t)=\eta^{*, \varepsilon}$ and $\eta^{*, \varepsilon}$, $b(t)=\beta_{1, \varepsilon}$ and $\beta_2^{\varepsilon}$, $k(t)=k^{\varepsilon}\triangleq k_0(\mu, \eta^{*, \varepsilon}, \beta_{1, \varepsilon})(t)$ and
 $k(t)=k_{\varepsilon}\triangleq k_0(\mu, \eta^{*, \varepsilon}, \beta_2^{\varepsilon})(t)$, respectively.

Since
\begin{align*}
&U_{\varepsilon}(t, r)\rightarrow \bar{V}_{\varepsilon}(t)\quad \mbox{in}~[0, T]\quad \mbox{as}~r\rightarrow+\infty,\\[5pt]
&V_{\varepsilon}(t, r)\rightarrow \underline{V}_{\varepsilon}(t)\quad \mbox{in}~[0, T]\quad \mbox{as}~r\rightarrow+\infty,
\end{align*}
there exists $R_0^*\triangleq R_0^*(\varepsilon)>2\hat{R}^*$ such that
\begin{align*}
&U_{\varepsilon}(t, r)>\bar{V}_{\varepsilon}(t)(1-\varepsilon)\quad\mbox{for}~(t, r)\in[0,T]\times[R_0^*, +\infty),\\[5pt]
&V_{\varepsilon}(t, r)>\underline{V}_{\varepsilon}(t)\quad\mbox{for}~(t, r)\in[0,T]\times[R_0^*, +\infty).
\end{align*}

We now define
\begin{align*}
&\bar{h}(t)=(1-\varepsilon)^{-2}\int_{0}^{1}k^{\varepsilon}(s)ds+R_0^*+\hat{h}(\hat{T}), \quad\mbox{for}~t\geq0,\\[5pt]
&\bar{u}(t, r)=(1-\varepsilon)^{-2}U_{\varepsilon}(t, \bar{h}(t)-r),\quad\mbox{for}~t\geq0, 0\leq r\leq \bar{h}(t),
\end{align*}
and
\begin{align*}
&\underline{h}(t)=(1-\varepsilon)^{-2}\int_{0}^{1}k_{\varepsilon}(s)ds+R_0^*+\hat{h}(0), \quad\mbox{for}~t\geq0,\\[5pt]
&\underline{u}(t, r)=(1-\varepsilon)^{-2}V_{\varepsilon}(t, \underline{h}(t)-r),\quad\mbox{for}~t\geq0, 0\leq r\leq \underline{h}(t).
\end{align*}
Direct calculations show that $(\underline{u}, \underline{h})$ and $(\bar{u}, \bar{h})$ are sub-supersolutions to $(1.1)$, respectively. Therefore, we have
\begin{align*}
\liminf_{t\rightarrow+\infty}\frac{h(t)}{t}
\geq\lim_{t\rightarrow+\infty}\frac{\underline{h}(t)}{t}=
(1-\varepsilon)^{-2}\frac{1}{T}\int_{0}^{T}k_{\varepsilon}(t)dt,\\[5pt]
\limsup_{t\rightarrow+\infty}\frac{h(t)}{t}
\leq\lim_{t\rightarrow+\infty}\frac{\bar{h}(t)}{t}=
(1-\varepsilon)^{-2}\frac{1}{T}\int_{0}^{T}k^{\varepsilon}(t)dt.
\end{align*}
Since $\varepsilon>0$ can be arbitrarily small, this implies our desired result.\quad$\Box$

The result below follows trivially from Theorem 6.1.

\noindent\textbf{Corollary 6.1.}
Assume that $h_{\infty}=+\infty$ and
\begin{align*}
\alpha(t, r)-\gamma(t, r)\rightarrow \breve{\eta}(t),\quad \beta(t, r)\rightarrow\breve{\beta}(t)\quad\mbox{uniformly for}~t\in[0, T]\quad\mbox{as}~r\rightarrow+\infty.
\end{align*}
Then
\begin{align*}
\lim_{t\rightarrow+\infty}\frac{h(t)}{t}=\frac{1}{T}\int_{0}^{T}k_0(\mu, \breve{\eta}, \breve{\beta})(t)dt.
\end{align*}

\end{document}